\RequirePackage{fix-cm}
\documentclass{svjour3}                     
\smartqed  

\usepackage[colorlinks,bookmarksopen,bookmarksnumbered,citecolor=red,urlcolor=red]{hyperref}

\usepackage{calrsfs}
\usepackage{subcaption}
\usepackage{graphicx}
\usepackage{epsfig}
\usepackage{multirow}
\usepackage{rotating}
\usepackage{bigstrut}
\usepackage{textcomp}
\usepackage[sumlimits]{amsmath}
\usepackage{amssymb}
\usepackage{mathrsfs}
\usepackage{appendix}
\numberwithin{equation}{section}
\numberwithin{figure}{section}
\usepackage{array}
\usepackage{booktabs}
\usepackage{graphicx}
\usepackage{epstopdf}
\usepackage[ansinew]{inputenc}
\allowdisplaybreaks[1]

\usepackage{framed,color}
\usepackage{color,soul}

\newcommand{\bzero}{{\mathbf 0}}

\newcommand{\eps}{{\varepsilon}}

\DeclareMathOperator{\tr}{tr} 
\DeclareMathOperator{\deviator}{dev} 
\DeclareMathOperator{\SymDev}{\mathrm{SymDev}} 

\newcommand{\bb}{\boldsymbol{b}}

\newcommand{\bu}{\boldsymbol{u}}

\newcommand{\balpha}{\boldsymbol{\alpha}}
\newcommand{\beps}{\boldsymbol{\varepsilon}}

\newcommand{\bsigma}{\boldsymbol{\sigma}}

\newcommand{\bH}{\boldsymbol{H}}
\newcommand{\bI}{\boldsymbol{I}}

\newcommand{\bK}{\boldsymbol{K}}

\newcommand{\e}{^\text{e}}
\newcommand{\ch}{^\text{ch}}

\newcommand{\trn}{^\text{tr}}

\newcommand{\etr}{\bbe\trn}

\newcommand{\y}{_{\textrm{y}_0}}

\newcommand{\dev}{_{\text{dev}}}

\newcommand{\odd}{\hspace{0.2mm}$\ddot{\mbox{o}}$}

\newcommand{\eldom}      {{\mathcal E}}
\newcommand{\real}       {\mathbb{R}}

\newcommand{\bbe}   {\mathbf{e}}

\newcommand{\bbs}   {\mathbf{s}}

\newcommand{\bbX}   {\mathbf{X}}

\newcommand{\temp}       {T}

\newcommand{\ratezeta}   {{\dot\zeta}}

\newcommand{\rateetr}    {{\dot{\bbe}}^{\mathrm{tr}}}

\newcommand{\Mf}        {M_{\textrm{f}}}

\newcommand{\epsL}   {\eps _\textrm{L}}
\newcommand{\satdom} {{\mathcal S}}

\newcommand{\bDTp}   {\beta\Delta T^{+}}

\newcommand{\Iepsl}  {{\cal I}_{\satdom}}

\def\wbox#1;#2;{\vbox{\hrule\hbox{\vrule height#1mm\kern#2mm\vrule height#1mm}\hrule}}

\renewcommand{\phi}{\varphi}

\def\F#1#2{{}_{#1}\kern-0.5ptF_{#2}}
\def\U#1#2{{}_{#1}\kern-0.5ptU_{#2}}

\def\odd#1{\kern-#1pt{}^{\textrm{o}}\,\,}

\newcommand\lshad{{[\kern-0.15em[}}
\newcommand\rshad{{]\kern-0.15em]}}

\newcommand{\vect}[1]{\bold{#1}}

\def\div{{\rm div}}

\def\x{{\bf x}}

\def\u{{\vect{ u}}}
\def\v{{\vect{ v}}}
\def\S{{\vect{ S}}}

\def\F{{\vect{ f}}}

\def\V{{\bf V}}
\def\W{{\bf W}}
\def\bD{{\bf D}}
\def\bu{{\bf u}}
\def\bI{{\bf I}}
\def\tr{\textrm{tr}}
\def\sgrad{\boldsymbol{\varepsilon}}

\def\eps{\boldsymbol{\varepsilon}}

\def\P{{\mathcal P}}  

\def\O{\Omega}

\def\Th{{\mathcal T}_h}
\def\E{E}

\def\VE{{\V_{h|\E}}}

\def\P{\mathsf P}

\definecolor{mygreen}{rgb}{0.2,0.8,0.2}

\def\u{{\bf u}}
\def\v{{\bf v}}

\def\s{{\bf s}}

\def\c{{\bf c}}

\def\sigmab{\boldsymbol{\sigma}}

\def\div{\textrm{div}}

\def\ah{a_h}

\def\dev{{\rm dev}}
\def\tr{{\rm tr}}

\def\Th{\Omega_h}

 \newcommand{\T}{^\text{T}}
\journalname{Computational Mechanics}
\begin{document}
\title{Arbitrary order 2D virtual elements for polygonal meshes: Part II, inelastic problem
\thanks{The first author gratefully acknowledges the partial financial support of the Italian Minister of University and Research, MIUR  (Program: �Consolidate the Foundations 2015�; Project: BIOART; Grant number (CUP): E82F16000850005). \\
The third author was partially supported by IMATI-CNR of Pavia, Italy. This support is gratefully acknowledged.}
}

\titlerunning{VEM for inelastic problems}        

\author{E. Artioli         \and
        L. Beir\~ao da Veiga \and
        C. Lovadina           \and
        E. Sacco   
}


\institute{
           E. Artioli \at
           Department of Civil Engineering and Computer Science\\
           University of Rome Tor Vergata\\
           Via del Politecnico 1, 00133 Rome, Italy\\	
           Tel.: +39-06-72597014\\
           Fax:  +39-06-72597005\\
           \email{artioli@ing.uniroma2.it}
           \and
           L. Beirao da Veiga \at
           Department of Mathematics and Applications\\
           University of Milan Bicocca\\
           Via Roberto Cozzi 55, 20125, Milan, Italy\\	
           Tel.: +39-02-64485773\\
           Fax:  +39-02-\\
           \email{lourenco.beirao@unimib.it}
           \and
           C. Lovadina \at
           Department of Mathematics\\
           University of Milano\\
           Via Cesare Saldini 50, 20133, Milano, Italy\\
           Tel.: +39-02-50316187\\
           Fax:  +39-02-50316090\\
           \email{carlo.lovadina@unimi.it}
           \and
           E. Sacco \at
           Civil Engineering and Mechanical Engineering\\
           University of Cassino and Southern Lazio\\
           Via Di Biasio 43, 03043 Cassino, Italy\\
           Tel.: +39-0776-2993659\\
           Fax:  +39-0776-2993392\\
           \email{sacco@unicas.it}
          }
                    
\date{Received: date / Accepted: date}

\maketitle

\begin{abstract}
The present paper is the second part of a twofold work, whose first part is reported in \cite{ABLS_part_I}, concerning a newly developed Virtual Element Method (VEM) for 2D continuum problems. 
The first part of the work proposed a study for linear elastic problem.
The aim of this part is to explore the features of the VEM formulation when material nonlinearity is considered, showing that the accuracy and easiness of implementation discovered in the analysis inherent to the first part of the work are still retained. 
Three different nonlinear constitutive laws are considered in the VEM formulation. In particular, the generalized viscoplastic model, the classical Mises plasticity with isotropic/kinematic hardening and a shape memory alloy (SMA) constitutive law are implemented.
The versatility with respect to all the considered nonlinear material constitutive laws is demonstrated through several numerical examples, also remarking that the proposed 2D VEM formulation can be straightforwardly implemented as in a standard nonlinear structural finite element method (FEM) framework.
\end{abstract}


\section{Introduction}
\label{s:intro}

The virtual element method has been introduced recently in \cite{volley,VEM-elasticity,Brezzi:Marini:plates,hitchhikers,projectors} as a generalization of the finite element method capable to deal with general polygonal/polyhedral meshes. The VEM approach has experienced an increasing interest in the recent literature, both from the theoretical (mathematical) viewpoint, and on the applicative (engineering) side. In an absolutely non-exhaustive way, in addition to the ones above we here limit to cite the few works \cite{BeiraoLovaMora,variable-primal,Benedetto-VEM-2,Berrone-VEM,BFM,Paulino-VEM,Gatica,VemSteklov,Helmo-PPR,wriggers,Vacca-1}. However, we note that VEM is not the only recent method that can make use of polytopal meshes, and we refer to \cite{Cangiani:Georgoulis:Houston2014,DiPietro-Ern-1,Droniou-gradient,VEM-topopt,Gillette-1,ST04,POLY37}, again without pretending to provide a complete picture of the available approaches on the topic. 
In the more specific framework of structural mechanics, VEM has been introduced in \cite{VEM-elasticity} for (possibly incompressible) two dimensional linear elasticity and general ``polynomial'' order, in \cite{Paulino-VEM} for three dimensional linear elasticity and lowest order, in \cite{BeiraoLovaMora} for general two dimensional elastic and inelastic problems under small deformations (lowest order), in \cite{wriggers} for contact problems and in \cite{Andersen-geo} for applications in geomechanics (again both contributions being for lowest order).

The present paper represents the continuation of the investigations started in \cite{ABLS_part_I}. In fact, the previous paper was devoted to a new VEM formulation for linear 2D elastic problems; the present paper is focused on the the extension of the proposed developed VEM formulation to problems with material nonlinearity. In particular, the aim of the paper is to suitably modify the VEM proposed in Part I to a general setting in which nonlinear inelastic constitutive behavior is taken into account, for arbitrary order of accuracy (or ``polynomial'' order).
Moreover, a numerical assessment of such VEM scheme is presented for three typical inelastic problems:
\begin{itemize}
	\item generalized Maxwell isotropic viscoelasticity; 
	\item classical von Mises plasticity with linear isotropic/kinematic hardening;
	\item shape memory alloy constitutive behavior modeled by means of a macroscopic phenomenological approach.
\end{itemize}
All the considered problems fit into a general framework, capable of modeling a wide class of inelastic effects, governed by a phenomenological constitutive law.
A similar approach using VEM technology has been initially presented in \cite{BeiraoLovaMora}, limited to the case of a low-order scheme. Here, we extend the higher order schemes presented in \cite{ABLS_part_I} for linear problems to material inelastic response. Analogously to \cite{BeiraoLovaMora}, one feature of the present approach is that the constitutive law algorithm can be independently embedded as a self-standing black-box, as in common nonlinear FEM codes. However, in addition to considering a general ``polynomial'' degree, the numerical tests presented in this paper generally differs from the ones provided in \cite{BeiraoLovaMora}, and some aspects concerning the computational behavior of the proposed VEM scheme are discussed.
Ultimately, the method is shown to be an appealing alternative for inelastic problem with respect to standard FEM.

An outline of the paper is as follows. 
In Sec. \ref{s:statement} the equilibrium problem for a 2D medium characterized by inelastic response, both in the continuous and in the VEM-discretized frameworks, is introduced. Purposely, the formulation is here kept quite general, in order to consider a wide gallery of constitutive models. Furthermore, we remark that implementation details of the proposed approach may be found in \cite{ABLS_part_I}, where a comprehensive discussion in the linear framework has been developed. The extension to solution of equilibrium equations in the non-linear case, using a Newton-Raphson strategy, follows standard steps, and it is not detailed in this paper for brevity.
Sec. \ref{s:constitutive} reviews three typical inelastic constitutive models, belonging to the general category recalled in Sec. \ref{ss:one}, which will be used in the numerical tests. Numerical results are given in Sec. \ref{s:numres}. Sec. \ref{s:conclusion} draws some conclusion and briefly present possible extensions of the schemes here proposed and studied.  
\section{Statement of the problem}
\label{s:statement}
%

\subsection{The continuous problem}
\label{ss:one}

In this section we present a quite general framework for inelastic problems in 2D, under the assumption of small strain and displacement. In the following the Voigt notation is adopted, so that stress and strain tensors are represented as $3-$component vectors, and the fourth-order constitutive tensor is represented as a $3\times3$ matrix.

Let $\O$ be a continuous body occupying a region of the two-dimensional space $\mathbb{R}^2$ in which the Cartesian coordinate system $(\textrm{O},x,y)$ is introduced.
The displacement field is denoted by the vector $\u(x,y)=\left \{ u \; v \right\}^T$ and the associated strain defined as:
\begin{equation}
\label{eq:epsilon}
\beps(\u) = \S \u \qquad \textrm{with} \;
\S = \left [ \begin{matrix} \partial_{x} & 0 \\ 0 & \partial_{y} \\ \partial_{y} & \partial_{x} \end{matrix} \right] \, .
\end{equation}
Above, the symbol $\partial_{(\bullet)}$ indicates the partial derivative operator with respect to the $(\bullet)$-coordinate.
An additive decomposition is considered for the total strain $\beps=\beps(\u)$, assuming the form:
\begin{equation}
\label{add-dec}
\beps=\beps^e+\beps^{\textrm{in}} ,
\end{equation} 
where $\beps^e$ is the elastic strain, and $\beps^{\textrm{in}}$ is the internal variable which represents the strain stemming from the inelastic effects. 

Setting the framework of generalized standard materials with convex free-energy and dissipation potential and following standard thermodynamic arguments, we consider a constitutive law for the body $\O$, such that the stress $\bsigma$ is given by the relationship:
\begin{equation}
\label{eq:sigma}
\bsigma = \sigmab(t,\x,\beps, \dot\beps, \beps^{\textrm{in}},\dot\beps^{\textrm{in}},\bH)  
\end{equation}
where $t$ is the time variable, $\x=\left\{x,y\right\} \T \in \O$ is the position vector, the vector $\bH$ contains all the history variables incorporated in the selected model, and a dot above a function stands, as usual, for the time derivative. 

\begin{remark}\label{genrule}
	We have defined the constitutive law \eqref{eq:sigma} in full generality, allowing for a rule depending on all the quantities involved in the description of a generalized standard material (but not temperature). However, we notice that in many interesting situations, $\sigmab$ depends only on $(\beps, \beps^{\textrm{in}},\bH)$ (see, for instance \cite{lubliner90,SouzaNeto}).
\end{remark}	

Rule \eqref{eq:sigma} is coupled, respectively, with an evolution law $\mathcal{L}$ for the inelastic strain, and an evolution law $\mathcal{M}$ for the history variables:
\begin{equation}
\label{evol-law}
\begin{array}{ll}
\dot \beps^{\textrm{in}}(t,\x) &\in {\cal L}(t,\x,\beps(t,\x),\dot\beps(t,\x), \beps^{\textrm{in}}(t,\x), \dot\beps^{\textrm{in}}(t,\x), \bH(t,\x) ) \\
\dot \bH(t,\x) &\in {\cal M}(t,\x,\beps(t,\x),\dot\beps(t,\x), \beps^{\textrm{in}}(t,\x), \dot\beps^{\textrm{in}}(t,\x), \bH(t,\x) ). 
\end{array} 
\end{equation}

\begin{remark}\label{mvalued}
Depending on the physical phenomenon under consideration, the evolution laws $\mathcal{L}$ $and \mathcal{M}$ may be described by either a standard single-valued correspondence, or a more general set-valued function. In particular, in several interesting situations $\mathcal{L}$ turns out to be associated with the sub-differential of a suitable {\em yield function}. For  more details we refer to \cite{Daya99} or \cite{JirasekBazant}, for instance.  
\end{remark}

\begin{remark}\label{rem:Mfunction}
In many constitutive models, for instance in the examples of Section \ref{s:constitutive}, the history variables  can be explicitly expressed as a function of the inelastic strains, say $\bH(t,\x)=\varphi(\beps^{\textrm{in}}(t,\x))$. In such cases, by taking the time derivative of this equation, it is possible to provide only the evolution law ${\cal L}$ (and the function $\varphi$, of course) to completely determine the model. In other words, the evolution law ${\cal M}$ can be computed by means of ${\cal L}$ and $\varphi$.     		
\end{remark}


We denote with $\bK_T$ the tangent matrix consistently computed from the constitutive law \eqref{eq:sigma}, i.e.:

\begin{equation}\label{tangent}
\begin{aligned}
\bK_T(t,\x,& \beps(t,\x),\dot\beps(t,\x), \beps^{\textrm{in}}(t,\x), \dot\beps^{\textrm{in}}(t,\x), \bH(t,\x))\\
&=\frac{\partial}{\partial\beps}\sigmab(t,\x,\beps(t,\x), \dot\beps(t,\x), \beps^{\textrm{in}}(t,\x), \dot\beps^{\textrm{in}}(t,\x), \bH(t,\x)) 
\end{aligned}
\end{equation}


The body $\O$ is subjected to distributed volume forces $\bb$. For simplicity, and without loss of generality, we assume that the displacements vanish on the whole boundary of $\O$. Since we consider a quasi-static problem, at each time instant the stresses and displacements must satisfy the equilibrium equations and boundary conditions, that read:

\begin{equation}\label{elast-prob-strong}
\left\{
\begin{aligned}
&  \div\, \sigmab + \bb = \mathbf{0}  \qquad & &\textrm{in} \ \Omega , \\
& \u = \mathbf{0} \qquad & &\textrm{on} \ \Gamma=  \partial\Omega . \\
\end{aligned}
\right.
\end{equation}

Now, let $\V$ denote the space of admissible displacements and $\W$ the space of its variations; both spaces will, in particular, satisfy the homogeneous Dirichlet boundary condition on $\Gamma$.
Assuming initial values $\beps^{\textrm{in}}_0(\x)$ and $\bH_0(\x)$ for the inelastic deformation and the history variables, respectively, a possible variational formulation of our inelastic problem can be written as:\\
\begin{equation}\label{inelast-prob}
\left\{
\begin{aligned}
& \textrm{For all } t \in (0,T], \textrm{ find } \u(t,\cdot) \in \V \textrm{ such that } \\
&\int_{\Omega} \sigmab(t,\x,\beps(t,\x), \dot\beps(t,\x), \beps^{\textrm{in}}(t,\x),  \dot\beps^{\textrm{in}}(t,\x), \bH(t,\x))^T \beps(\v(\x)) d\x\\ 
& \hskip2cm =  \int_{\Omega} \bb(t,\x)^T\v(\x) d\x
\quad \forall \v \in \W \\
& \beps^{\textrm{in}}(0,\x)=\beps^{\textrm{in}}_0(\x)\qquad \forall \x\in\O  \\
& \bH(0,\x)=\bH_0(\x)\qquad \forall \x\in\O ,
\end{aligned}
\right.
\end{equation}
where the displacements and history variables are sufficiently regular in time and must satisfy the evolution laws \eqref{evol-law}. 


\subsection{The virtual element formulation}
\label{s:abs}

We now describe the virtual element method when applied to the problem class described in Sec. \ref{ss:one}. We will closely follow the notations and the framework of \cite{ABLS_part_I}, where the virtual element philosophy for the easier problem of linear elasticity, has been extensively detailed.

We start by presenting the discrete (virtual) space of admissible displacements $\V_h$, which is the same of \cite{VEM-elasticity}, see also \cite{ABLS_part_I}.
As discrete space for the displacement variations, we choose $\W_h := \V_h$.
Let $\Th$ be a {\it simple polygonal mesh} on $\O$, i.e. any decomposition of $\O$ into non-overlapping polygons $\E$ with straight edges. 
The symbol $m$ represents the number of edges of a polygon $E$, and the typical edge of the polygon $E$ is indicated by $e$, (i.e. $e \in \partial E$).
The space $\V_h$ will be defined element-wise, by introducing local spaces $\VE$ and the associated local degrees of freedom, as in standard Finite Element (FE) analysis. On the other hand, differently from standard FE, the definition of the local spaces $\VE$ is not fully explicit. 

Let $k$ be a positive integer, representing the ``degree of accuracy'' of the method. Then, given an element $\E\in\Th$, we define:
\begin{equation}
\begin{aligned}
\VE=\big\{ & \v_h \in [H^1(\E) \cap C^0(\E)]^2 \ : \  \Delta\v_h   \in [\P_{k-2}(\E)]^2 , \\
& \v_h|_e \in [\P_k(e)]^2   \quad \forall e \in \partial\E \big\} ,
\end{aligned}
\end{equation}
where, for any subset $F \subseteq \O$, $\P_k(F)$ is the space of polynomials on $F$ of degree $\le k$, with the agreement that $\P_{\!-1} = \{ 0 \}$.

The space $\VE$ is made of vector valued functions $\v_h$ such that:
\begin{itemize}
\item $\v_h$ is a polynomial of degree $\le k$ on each edge $e$ of $\E$, i.e. $\v_h\in[\P_k(\E)]^2$;
\item $\v_h$  is globally continuous on $\partial\E$;
\item the laplacian $\Delta\v_h$ is a polynomial of degree $\le k-2$ in $\E$.
\end{itemize}

%
For the dimension of the space $\VE$, it holds:
\begin{equation}
\label{dofsnumber}
\textrm{dim}(\VE) = 2mk + k (k-1) .
\end{equation}


As in standard FE methods, the global space $\V_h \subseteq \V$ is built by assembling the local spaces $\VE$ as usual:
\begin{equation}
\V_h = \{ \v \in \V \ : \ \v |_{\E} \in \VE \; \forall \E\in\Th \} .
\end{equation}

We now introduce a {\it projection} operator $\Pi$ 
\begin{align}
\label{uVuP}
\Pi :  \;\;\;\; & \VE \longrightarrow \P_{k-1}(\E)^{2\times 2}_{\textrm{sym}} \\
\nonumber
& \quad \v_h \; \mapsto \quad \Pi(\v_h) 
\end{align}
to approximate the strain field $\beps^e(\v_h)$ induced by the virtual displacement $ \v_h $. 
Hence, for $\v_h \in \VE$, $ \Pi(\v_h) \in \P_{k-1}(\E)^{2\times 2}_{\textrm{sym}}$ is defined by:
\begin{equation}
\int_{\E} \Pi(\v_h)^T  \sgrad^{P}  = 
\int_{\E}      \eps (\v_h)^T  \sgrad^{P}  ,\quad \forall \sgrad^{P} \in \P_{k-1}(\E)^{2\times 2}_{\textrm{sym}} \,.
\label{eq:def-proj}
\end{equation}
This operator represents the best approximation of the strains (in the square integral norm) in the space of piecewise polynomials of degree $k-1$. We refer the reader to \cite{ABLS_part_I} for details on the construction and implementation of the operator $\Pi$.

In order to solve the constitutive evolution equation detailed in the next section, a Euler time integration is performed. To this end, we introduce a sub-division of the time interval $[0,T]$ into smaller intervals $[t_{n},t_{n+1}]$ for $n=0,1,...,N-1$, such that the time step is defined by $\Delta t= t_{n+1},t_{n}$. Correspondingly, partial loadings evaluated at $t_{n}$ are denoted as $\bb^n = (n/N) \bb$ for all $n=0,1,...,N-1$.

We assume, as in standard engineering procedures, a constitutive algorithm that is an approximation of the constitutive and evolution laws \eqref{eq:sigma}, \eqref{evol-law}. In Finite Element analysis, this pointwise algorithm can be coded independently from the global FE construction and can be regarded as a ``black-box'' procedure that is applied at every Gauss point and at every iteration step.
In the present Virtual Element method, we want to keep the same approach; in other words, our scheme will be compatible with any black-box constitutive algorithm that falls in the general setting below and that can be imported from other independent sources.

Let $\widehat\sigmab$ represent the constitutive algorithm, in the framework of the strain driven procedure based on a  backward-Euler approach. 
Hence, given:
\begin{enumerate}
\item a value for the strain $\beps_h^{n}(\x)$ at time $t_{n}$, 
\item a value for the inelastic strain $(\beps^{\textrm{in}})^{n}(\x)$ at time $t_{n}$, 
\item a value $\bH^{n}(\x)= \bH(t_{n},\x)$ for the internal variables at time $t_{n}$, 
\item a tentative value for the strain $\beps_h^{n+1}(\x)$ at time $t_{n+1}$, 
\end{enumerate}
the algorithm computes the stresses  at time $t_{n+1}$, and updates the inelastic strains and the history variables (i.e. it returns also $(\beps^{\textrm{in}})^{n+1}(\x)$ and $\bH^{n+1}(\x)$).
We thus write the computed stress at time $t_{n+1}$ and spatial location $\x$ as
\begin{equation}\label{blackbox}
\widehat\sigmab^{n+1}(\x) = \widehat\sigmab(t_{n+1}, \x, \beps_h^{n}(\x), (\beps^{\textrm{in}})^{n}(\x),  \bH^{n}(\x), \beps_h^{n+1}(\x)) .
\end{equation}
pointing out that the functional dependence of the updated stress tensor $\widehat\sigmab^{n+1}(\x)$ on quantities evaluated at time $t_{n}$ (and not only at time $t_{n+1}$) is to be viewed in algorithmic sense.

For a given element $E\in \Th$, we now select a suitable set of $np=np(E)$ points $\{ \x_{i,E}  \} \subset E$, for $i = 1,\cdots, np$. These points may be seen as the VEM analogous to the Gauss points for developing the numerical quadrature in standard Finite Elements. 
We then denote with $\bH^k_E$ the vector collecting the values $\{\bH(t_k, \x_{i,E})\}_{i=1}^{np}$, for $k=1,2,...,N$. 
Similarly, we set $\bH^k$ as the vector collecting all the vectors $\bH^k_E$, with $E\in\Th$

The Virtual Element scheme reads, for $n=0,2,...,N-1$:
\begin{equation}\label{inelast-virt-prob}
\left\{
\begin{aligned}
& \textrm{Find } \u_h^{n+1} \in \V_h \textrm{ (and the updated } \bH^{n+1}) \textrm{ such that } \\
& \ah(\u_h^{n}, \bH^{n}; \u_h^{n+1}, \v_h)  = <\bb^{n+1} , \v_h>_h \quad \forall \v_h \in \V_h .
\end{aligned}
\right.
\end{equation}

Above, the form $\ah(\u_h^{n}, \bH^{n}; \u_h^{n+1}, \v_h)$ is the sum of local contributions:

\begin{equation}\label{ah-form}
\ah(\u_h^{n}, \bH^{n}; \u_h^{n+1}, \v_h) = \sum_{E\in\Omega_h} \ah^E(\u_h^{n}, \bH^{n}_E; \u_h^{n+1}, \v_h) .
\end{equation}

To emphasize the dependence on the displacement field $\u_h^{n+1}$, we set (cf. \eqref{blackbox}):

\begin{equation}\label{sigmai}
\begin{aligned}
 \widehat\sigmab^{n+1}(\x_{i,E}, \u_h^{n+1}) &:=  \\
 &\widehat\sigmab\left( t_{n+1}, \x_{i,E}, \beps_h^{n}(\x_{i,E}), (\beps^{\textrm{in}})^{n}(\x_{i,E}),\bH^{n}(\x_{i,E}), \beps_h^{n+1}(\x_{i,E}) \right) ,
\end{aligned}
\end{equation}
where, given $\Pi$ defined by \eqref{eq:def-proj}, we have:

\begin{enumerate}
\item $\beps_h^{n}(\x_{i,E}):= \Pi(\u_h^{n})(\x_{i,E})$ is the computed total strain at $(t_{n},\x_{i,E})$;
\item $(\beps^{\textrm{in}})^{n}(\x_{i,E})$
is the computed inelastic strain at $(t_{n},\x_{i,E})$;
\item $\bH^{n}(\x_{i,E})$
are the computed history variables at $(t_{n},\x_{i,E})$;
\item $\beps_h^{n+1}(\x_{i,E}):= \Pi(\u_h^{n+1})(\x_{i,E})$ is the unknown total strain at $(t_{n+1},\x_{i,E})$.
\end{enumerate}

The local form is thus given by

\begin{equation}\label{ahE-form}
\begin{aligned}
 \ah^E(\u_h^{n},  \bH^{n}_E; & \u_h^{n+1}, \v_h) : = \\
& \sum_{i=1}^{np} \omega_{i,E}\, \widehat\sigmab^{n+1}(\x_{i,E},\u_h^{n+1})^T \Pi(\v_h^{n+1})(\x_{i,E})  + S^E \left( \u_h , \v_h \right) ,
\end{aligned}
\end{equation}
%
%
%
where $\{ \omega_{i,E} \}_{i=1}^{np}$ is a suitable set of weights, and $S^E \left( \u_h , \v_h \right)=\alpha(E)\,s^E \left( \u_h , \v_h \right)$ is a stabilization term, similar to the one involved in the linear case, see \cite{ABLS_part_I}. However, as investigated in \cite{BeiraoLovaMora}, in the present inelastic quasi-static setting the parameter $\alpha(E)$ needs to be differently chosen, to improve the robustness of the method. More precisely, we here select the mean value of the trace of the tangent matrix computed at time $t_n$, cf. \eqref{tangent}:

\begin{equation}\label{alpha-inel}
\alpha(E) = \frac{ \tr\left( \bK_T(t_{n},\c_E,\Pi(\u_h^{n}), \bH(t_{n},\c_E))\right)}{m} ,
\end{equation}
where $\c_E$ is the barycenter of $E$. 

Finally, the computation of the loading term $<\bb^{n+1} , \v_h>_h$ in \eqref{inelast-virt-prob} follows exactly the guidelines of \cite{ABLS_part_I,VEM-elasticity}.

\begin{remark}\label{approxint}
We remark that the local form in \eqref{ahE-form} is to be intended as the approximation of the energy integral
$\displaystyle{\int_E \bsigma^T\beps}$ over the polygon $E$. Consequently, the form in \eqref{ah-form} represents an approximation of the global internal energy
$\displaystyle{\int_\O \bsigma^T\beps}$ over the whole domain $\Omega$.
\end{remark}

\begin{remark}\label{adapt-gp}
The points $\{ \x_{i,E}  \}_{i=1}^{np}$ and the weights $\{ \omega_{i,E} \}_{i=1}^{np}$ 
must be chosen to make the associated integration rule exact for polynomials of degree up to $2(k-1)$, as it happens for standard triangular Finite Elements of degree $k$. Since we are here treating general polygons, such rules can be built, for instance, either by using a coarse sub-triangulation (that in the convex case is very easy to build), or by adopting more specific approaches, see \cite{MS11}. 
Note moreover that, to make the presentation as easy as possible, the selection of points $\{ \x_{i,E}  \} \subset E$, for $i = 1,\cdots, np(E)$ is here independent of the time variable. However, we remark that, in practice, the points $\{ \x_{i,E}  \}$ might
suitably vary at each time instant, following a sort of an ``adaptive'' strategy.  
\end{remark}
\section{Constitutive models} \label{s:constitutive}

In this section, a set of phenomenological nonlinear constitutive models are briefly presented, in the so called {\it energetic format} (see \cite{Mielke}, for instance), in order to give discussion and applications a unified layout. 
The first model is the generalized Maxwell viscoelastic model \cite{Zienkiewicz_Taylor_Zhu13}, then the classical von Mises plasticity model with linear isotropic and kinematic strain hardening is illustrated \cite{simo_computational_1998}. Finally, the shape memory alloy model proposed in \cite{Souza98} and, then, modified in \cite{auricchio02a,Evangelista2009} is presented. The introduced models are chosen to verify the effectiveness of the VEM methodology in reproducing classical nonlinear effects such as viscoelasticity plasticity, and shape memory of structural elements, and to prove superior behavior in such instances with respect to standard displacement-based finite element schemes.

\subsection{Generalized Maxwell isotropic viscoelastic constitutive model}
\label{ss:viscoelmodel}
The considered constitutive model is comprised of a linear elastic element in parallel with $M$ spring-dashpot linear elements, leading to a Helmholtz internal energy density of the following kind:
\begin{eqnarray}
\label{eq:psi_visc}
\psi\e(\beps, \mathbf{q}{^{(m)}}) &=& \frac{1}{2} \beps{^{(0)}} \bD^{(0)} \beps{^{(0)}} + \frac{1}{2} \sum_{m=1}^M \mathbf{q}{^{(m)}} \bD^{(m)} \mathbf{q}{^{(m)}}
 \end{eqnarray}
where $\beps{^{(0)}} = \beps$ and $\bD^{(0)}$ are the strain and linear elasticity matrix associated with the single elastic element; the terms $\mathbf{q}{^{(m)}}$ and $\bD^{(m)}$, $m=1,...M$, are the {\it partial strains} and elasticities in the dissipative spring-dashpots elements \cite{Zienkiewicz_Taylor_Zhu13}.

Applying standard continuum thermodynamics, the constitutive equation is derived \cite{Zienkiewicz_Taylor_Zhu13}:
\begin{eqnarray}
\label{eq:sigma-visco}
\mathbf{\bsigma}(t) & = & \bD^{(0)} \beps{^{(0)}}(t) + \sum_{m=1}^M \bD^{(m)} \mathbf{q}{^{(m)}}(t)   
\end{eqnarray}
In the above equation, each partial strain $\mathbf{q}{^{(m)}}(t) $ evolves according to:
\begin{equation}
\label{eq:visco-evol}
\dot{\mathbf{q}}{^{(m)}} + \frac{1}{\lambda^{(m)}} \mathbf{q} {^{(m)}}  = \dot{\beps}   
\end{equation}
where the terms $\lambda^{(m)}$ are coefficients of relaxation. In an integral form, the stress-strain behavior may be described through a convolution form as: 
\begin{equation}
\mathbf{\bsigma}(t) = \bD(t) \beps(0) + \int_{0}^t{\bD(t-\tau) \, \dot{\beps} \, d\tau} \,.\label{eq:st}
\end{equation}
where components of $\bD(t)$ are {\it relaxation moduli} functions.

Assuming isotropic material behavior, and considering a purely deviatoric inelastic response, the above relations simplify according to:  
\begin{eqnarray}
\label{eq:s-visco}
\mathbf{\s}(t) & = & 2 G(t) \mathbf{e}   
\end{eqnarray}
with $\mathbf{\s} = \dev \bsigma$, and $\mathbf{e} = \dev \beps$, and with $G(t)$ defined as the {\it shear modulus relaxation function}. 

In integral form the constitutive behavior is described as:
\begin{equation}
\mathbf{\s}(t) = \int_{-\infty}^t{2 G(t-\tau) \, \dot{\mathbf{e}} \, d\tau}
\end{equation}
The integral equation form may be defined as a generalized Maxwell model by assuming the {\it shear modulus relaxation function} in Prony series form \cite{Zienkiewicz_Taylor_Zhu13}:
\begin{equation}
\label{eq:Gt}
G(t) = G\,\left(\mu_0  + \sum_{i=1}^M  \mu_i \, \exp (- t / \lambda_i) \right) \,.
\end{equation}
\begin{remark}
\label{rem:vsel}
It is noted that the present constitutive model falls in the general framework outlined at the beginning of Sec. \ref{ss:one}. In particular, inelastic strains are given by the collection $\beps^{in}=\{\beps^{in,m}\} = \{\mathbf{q}{^{(m)}}\}$, $m = 1,...,M$ i.e. the partial strain tensors, while $\bH = \emptyset$ as no history variables are considered. Moreover, in every spring-dashpot element, the total strain is additively split in elastic and viscous parts cf. \eqref{add-dec}. The evolution law for the partial strains (cf. \eqref{eq:visco-evol}) is represented by a standard single-valued correspondence of viscous type, cf. \eqref{eq:visco-evol}:
\begin{equation}
\label{eq:evol_spec_vsel}
{\cal L}{^{(m)}} =
{\cal L}{^{(m)}(\dot\beps,\beps^{in,m}) = \dot{\beps}   - \frac{1}{\lambda^{(m)}} \beps^{in,m}} .
\end{equation}
Accordingly, the evolution law for the inelastic strains are given by the collection, cf. \eqref{evol-law}: 
\begin{equation}\label{eq:evol_spec_vsel-tot}
\begin{aligned}
&{\cal L}={\cal L}(\dot\beps,\beps^{in})=
\{ {\cal L}^{(m)}(\dot\beps,\beps^{in,m}) \}_{m=1}^{M}\\
&{\mathcal M}=\emptyset .
\end{aligned}
\end{equation}

\end{remark}

\subsection{Plasticity model}
\label{ss:plastmodel}
The von Mises plasticity model with combined linear isotropic/kinematic hardening is considered \cite{auricchio99b}.  

The strain is split into the deviatoric, $\mathbf{e}$, and volumetric (spherical), $\theta$, parts resulting:
\begin{equation}
\boldsymbol{\varepsilon} = \mathbf{e} + \frac{1}{2} \theta \bI \,,
\end{equation}
where $\mathbf{e} = \dev \beps$, $\theta=\tr \beps$.
Both the deviatoric and spherical strains are decomposed in the elastic and plastic parts:
\begin{eqnarray}
\mathbf{e} & = & \mathbf{e}\e +\mathbf{e}^{\mathrm{p}}        \label{eq:plastic-add-dec}\\
\theta         & = & \theta^e                   \label{eq:theta} \,.
\end{eqnarray}
indicating a purely isochoric plastic flow.

For an isotropic material, the Helmholtz free energy density assumes the form:
\begin{equation}
\psi = \psi\e + \psi\trn
\end{equation}
where:
\begin{eqnarray}
\label{eq:psi_ee}
\psi\e(\beps\e) &=& \frac{1}{2}K(\theta^e)^2 +
 G \| \mathbf{e}\e \|^2 \\
 \psi\trn (\mathbf{e}^{\mathrm{p}}) &=& \frac{1}{2} H^{\mathrm{kin}}\| \mathbf{e}^{\mathrm{p}} \|^2
 \end{eqnarray}
where $K$ and $G$ are, respectively, the bulk and shear elastic moduli, and $H^{\mathrm{kin}}$ is the linear kinematic hardening parameter.

By standard thermodynamic arguments, the constitutive equations are derived as: 
\begin{eqnarray}
\label{eq:sig}
&\boldsymbol{\s}  &\in \partial_{\mathbf{e}\e} \psi \\
\label{eq:relst}
&\bbX  &\in - \partial_{\mathbf{e}^{\mathrm{p}}}\psi\,,
\end{eqnarray}
Here, the symbol $\partial$ represents the subdifferential operator in the sense of Convex Analysis. The {\it relative stress} $\bbX$ in Eq. \eqref{eq:relst} is usually rewritten as:
\begin{equation}
\label{eq:bbX_split}
\bbX = \bbs - \balpha
\end{equation}
where $\bbs$ is the stress deviator, and $\balpha$ is the back stress tensor, which are given, respectively, by:
\begin{eqnarray}
\label{eq:s}
\mathbf{s}&=&2G\,\left(\mathbf{e} - \mathbf{e}^{\mathrm{p}} \right)   \\
\label{eq:backst}
\balpha &=& H^{\mathrm{kin}} \mathbf{e}^{\mathrm{p}} 
\end{eqnarray}

The activation of the plastic flow is governed by the von-Mises yield function  expressed in terms of the relative stress:
\begin{eqnarray}
\label{eq:f-mises}
f\left(\mathbf{X},\bar{e}^\mathrm{p}\right) &=& \| \bbX \| - \frac{\sqrt{2}}{2}\sigma_{\mathrm{y}}
\end{eqnarray}
which defines the {\it elastic domain} as the set $\eldom = \left \{\bbX \in \SymDev : f(\bbX) \le 0 \right \}$. The 
function $\sigma_{\mathrm{y}} = \sigma\y + H^{\mathrm{i}} \bar{e}^\mathrm{p}$ is the uniaxial yield stress, depending on the initial 
 yield stress $\sigma\y$, on the isotropic hardening parameter $H^{\mathrm{i}}$, and on the accumulated plastic strain
 
\begin{equation}\label{eq:accumul-def}
  \bar{e}^\mathrm{p} = \int_{0}^{t}{\| \dot{\mathbf{e}}^\mathrm{p} \| }d \tau .
\end{equation}  
The evolution law for the plastic strain tensor is {\it associated} to the yield function inasmuch it results:
\begin{equation}
\label{eq:plastic_flow}
\dot{\mathbf{e}}^\mathrm{p}  = 
 \dot{\zeta}\nabla_{\bbX} f
\end{equation}
from which it follows

\begin{equation}\label{eq:accumul-zeta}
\dot{\bar{e}}^\mathrm{p} = \frac{\sqrt 2}{ 2}\dot{\zeta} .
\end{equation} 
The evolution is complemented by the Kuhn-Tucker optimality conditions:
\begin{equation}
\dot{\zeta}\geq0\qquad\qquad \dot{\zeta}\,f=0 \,.\label{eq:KT-pla}
\end{equation}
for the plastic rate parameter $\dot{\zeta}$.
\begin{remark}
\label{rem:plas}
It is noted that the present constitutive model falls in the general framework outlined at the beginning of Sec. \ref{ss:one}. In particular, the inelastic strain is the plastic strain, i.e. $\beps^{\textrm{in}} = \mathbf{e}^{\mathrm{p}}$; the history variable is simply $\bH = \balpha $, and the total strain is still additively split into elastic and plastic parts, cf. \eqref{add-dec}. Using relations \eqref{eq:bbX_split}, \eqref{eq:accumul-zeta} and \eqref{eq:backst}, the evolution law in this case is given by (cf. \eqref{evol-law} and \eqref{eq:plastic_flow}):

\begin{equation}
\label{eq:evol_spec_plas}
\begin{aligned}
&{\cal L} =
{\cal L}(\beps, \beps^{\textrm{in}}, \dot\beps^{\textrm{in}}, \bH ) = \sqrt 2\ \| \dot\beps^{\textrm{in}} \| \
\frac{2G\,\dev \left(\beps - \beps^{\textrm{in}} \right) - \bH}{\| 2G\,\dev \left(\beps - \beps^{\textrm{in}} \right) - \bH \|}\\
&{\cal M} =
{\cal M}(\dot\beps^{\textrm{in}}) = H^{\mathrm{kin}}\, \dot\beps^{\textrm{in}} .
\end{aligned}
\end{equation}

We also notice that, due to Eq. \eqref{eq:backst} and Remark \ref{rem:Mfunction}, we could have provided the only evolution law ${\cal L}$ to describe the model.

\end{remark}

\subsection{Shape memory alloy constitutive model}
\label{ss:SMA}
The local thermodynamic state of the material is defined by the infinitesimal strain $\beps$,
the absolute temperature $\temp$, and by the symmetric {\it transformation} strain tensor $\etr$,
assuming the strain additive decomposition:
\begin{equation}
\label{eq:strain_decomp}
\beps = \beps\e + \etr
\end{equation}
into elastic strain, $\beps\e$, and transformation strain.
The quantity $\etr$ is the inelastic strain associated
with the phase transformation, assumed to be traceless indicating phase transition to be isochoric \cite{orgeas96}. The transformation strain is constrained to belong to the {\it saturation} domain $\satdom = \left \{\etr \in \SymDev : c(\etr) \le 0 \right \}$, where $c(\etr) = \|\etr\|^2/\epsL^2 - 1$, and $\epsL$ is a material parameter related to the maximum
transformation strain reached at the end of the forward isothermal
transformation during a uniaxial test. Given its tensorial character, $\etr$ is capable of representing the reorientation of the product
phase in the saturated condition \cite{auricchio02a}.

The Helmholtz free energy density $\psi$ is here taken as a strictly convex potential depending on the local thermodynamic state of the material:
\begin{equation}
\label{eq:free_energy}
\psi( \beps\e, \etr, \temp) = \psi\e(\beps\e) + \psi\ch(\etr,\temp) + \psi\trn(\etr) \,,
\end{equation}
under the constraint $\etr \in \satdom$. Here:
\begin{itemize}
\item
$\psi\e$ is the elastic strain energy, which, assuming linear isotropic elastic behavior, is given by: 
\begin{equation}
\label{eq:psi_e}
\psi\e(\beps\e) =\frac{1}{2}K(\tr\beps\e)^2 +
 G \| \deviator \beps\e \|^2
 \end{equation}
with  $K$ the bulk modulus and $G$ the shear modulus;
\item
$\psi\ch$ is the chemical energy, associated with the thermally-induced martensitic transformation:
\begin{equation}
\label{eq:psi_ch}
\psi\ch(\etr,\temp) = \bDTp \|\etr\|
\end{equation}
with $\beta$ a material parameter related to the dependence of the critical stress
on the temperature, and $\Delta T^{+} = \left<\temp - \Mf\right>$, being $\Mf$ the temperature corresponding to the end of the forward transformation, and $\left<\bullet\right>$ the positive part of the argument;
\item
$\psi\trn$ is the transformation strain energy, associated with transformation-induced strain hardening:
\begin{equation}
\label{eq:psi_tr}
\psi\trn(\etr) = \frac{1}{2}h {\|\etr\|}^2
\end{equation}
with $h$ a material parameter defining the slope of the linear stress -  transformation strain relation in the uniaxial case.

\end{itemize}
As a consequence of the principle of maximum inelastic dissipation \cite{Daya99}, the thermodynamic equilibrium state is expressed in terms of the quantities thermodynamically conjugate to the arguments $(\beps\e,\etr,\temp)$. By definition:
\begin{equation}
\label{eq:bbX}
\begin{aligned}
\bsigma & =  \partial_{\beps\e}\psi \,, \\
\bbX & = - \partial_{\etr}\psi - \partial_{\etr}\Iepsl\,, \\
\eta & = - \partial_{\temp}\psi \,.
\end{aligned}
\end{equation}
where $\bsigma$ is the Cauchy stress, $\bbX$ is the symmetric traceless thermodynamic stress, and $\eta$ is the entropy density.
Eq. \eqref{eq:bbX}$_2$ is usually rewritten as:
\begin{equation}
\label{eq:bbX_relative}
\bbX = \bbs - \balpha
\end{equation}
where $\bbs$ is the stress deviator, and $\balpha$ is the back stress tensor, given by:
\begin{equation}
\label{eq:alpha}
\balpha = \bDTp \partial_{\etr} \| \etr \| + h \, \etr +
\partial_{\etr} \Iepsl(\etr)
\end{equation}
Moreover, the indicator function $\Iepsl(\etr)$ of the saturation domain is introduced:
\begin{equation}
\label{eq:indicator}
\left\{
\begin{array}{ll}
\Iepsl(\etr) = 0   & \mbox{ if } \; c(\etr) \le 0
\vspace*{2mm} \\
\Iepsl(\etr) = +\infty & \mbox{ otherwise }
\end{array}
\right.
\end{equation}
whose subdifferential results:
\begin{equation}
\label{eq:indicator_subdiff}
\partial_{\etr} \Iepsl(\etr) =
\left\{
\begin{array}{ll}
\bzero & \mbox{ if } {c(\etr) < 0}
\vspace*{2mm} \\
\gamma \, \partial_{\etr} c(\etr) & \mbox{ if } {c(\etr) = 0}
\vspace*{2mm} \\
\emptyset  & \mbox{ if } {c(\etr)} > 0
\end{array}
\right.
\end{equation}
being $\gamma\in\real_0^{+}$ the Kuhn-Tucker parameter associated with the thermodynamic reaction explicated by the saturation constraint.

The phase transformation mechanism is governed by a flow law for the transformation strain, assigned in terms of a {\it transformation function} $f(\bbX)$ which defines the set of admissible thermodynamic stresses as the nonempty, closed {\it elastic domain} in deviatoric stress space\footnote{In the following, the space of symmetric traceless second-order tensors is denoted by $\SymDev$.}:
\begin{equation}
\label{eq:eldom}
\eldom = \left \{\bbX \in \SymDev : f(\bbX) \le 0 \right \}
\end{equation}
The function $f(\bbX)$ is typically expressed in the guise of a plasticity yield function, deviatoric isotropic, and represented in this context in von-Mises form \cite{Artioli_Bisegna_IJNME_2016a}: 
\begin{equation}
\label{eq:transf_func}
f(\bbX) = || \bbX || - \sqrt{\frac{2}{3}}\sigma\y
\end{equation}
Accordingly, activation of inelastic flow is ruled as follows:
\begin{eqnarray}
\left \{
\begin{array}{ll}
\mathrm{if} \; f(\bbX)<0, & \rateetr = \bzero \\
\mathrm{if} \; f(\bbX)=0, & \rateetr = \bzero, \mathrm{OR}\; \rateetr \neq \bzero.
\end{array}
\right.
\end{eqnarray}
In particular, the flow law is obtained by postulating the principle of maximum inelastic (or transformation) work rate \cite{Daya99,Hackl2008b}, implying convexity of the elastic domain, and leading to an associated flow rule in the form:
\begin{equation}
\label{eq:sma_flow}
\rateetr = \ratezeta \nabla f(\bbX)
\end{equation}
being $\bbX$ the admissible thermodynamic stress at equilibrium, with the Kuhn-Tucker conditions for the inelastic rate parameter $\ratezeta$:
\begin{equation}
\label{KTcond}
\ratezeta \ge 0 \;, \;\; \ratezeta f = 0 \;.
\end{equation}

\begin{remark}
\label{rem:smas}
It is noted that the present constitutive model falls in the general framework outlined at the beginning of Sec. \ref{ss:one}. In particular, the inelastic strain is the transformation strain, i.e. $\beps^{\textrm{in}} = \mathbf{e}^{\mathrm{tr}}$; the history variable is simply $\bH = \balpha $, and the total strain is still additively split into elastic and transformation parts, cf. \eqref{eq:strain_decomp}. Due to Eq. \eqref{eq:alpha}, the model description is determined by providing only the evolution law ${\cal L}$ for the inelastic strains, cf. Remark \ref{rem:Mfunction}.
Using relations \eqref{eq:bbX_relative} and \eqref{eq:transf_func}, such an evolution law in this case is given by (cf. \eqref{evol-law} and \eqref{eq:sma_flow}):

\begin{equation}
\label{eq:evol_spec_smas}
{\cal L} =
{\cal L}(\beps, \beps^{\textrm{in}}, \dot\beps^{\textrm{in}}, \bH ) = \sqrt 2\ \| \dot\beps^{\textrm{in}} \| \
\frac{2G\,\dev \left(\beps - \beps^{\textrm{in}} \right) - \bH}{\| 2G\,\dev \left(\beps - \beps^{\textrm{in}} \right) - \bH \|}.
\end{equation}
\end{remark}
\section{Numerical results}
\label{s:numres}
The present section is devoted to validation of the proposed VEM formulation in conjunction with the inelastic constitutive models previously examined (see Sec. \ref{s:constitutive}). A set of classical benchmarks is presented in order to assess accuracy of the proposed approach in comparison with standard Lagrangian displacement finite element schemes and hence to show the good properties of the VEM in regard to robustness and capability of efficient treatment of material nonlinearity. Noteworthy, the main tools of nonlinear finite element analysis of continua and structures such as Newton's method for solving equilibrium equations are still retained in the present context, highlighting the versatility of the VEM methodology.

\subsection{On the integration rule for polygons}
\label{ss:intrule}
This introductory section precedes the numerical test campaign inasmuch it gives some hints on the actual implementation of the method outlined in Sec. \ref{s:abs}. In particular, attention is centered on the computation of the residual equation which, in turn, entails computation of stress and material tangent stiffness at quadrature points over a polygon. A key point in such a procedure is computation of area integrals over polygonal domains. The problem obviously refers to cases $k \ge 2$ \footnote{Note that, in the case $k=1$, a single Gauss point for the whole polygon (for instance at the centroid of the element) is sufficient, see \cite{BeiraoLovaMora}.}. In the present context, given the convexity of all polygons in any mesh, a mere sub-triangulation of each polygon in $m$ triangles is adopted, by choosing the centroid as the common vertex shared by all triangles and tracing rays from such point to the vertexes. By doing so, any area integral can be computed through summation of the integral contributions stemming from each sub-triangle, computed, for instance, adopting Gaussian quadrature for triangles \cite{Zienkiewicz_Taylor_Zhu13} and leading to a total of (at most) $ms$ integration points on each polygon, being $s$ the number of (interior) Gauss points per sub-triangle. Of course, this integration strategy presumes to process the entire mesh element-wise in the integral computation loop. 

A different convenient strategy is to use a quadrature formula for sub-triangles with integration points also on the boundary of the element and not only in the interior. For example, in the case $k=2$, a convenient strategy is to use a quadrature formula for sub-triangles with a single integration point on each edge (at the midpoint), amounting to a total of $2m$ integration points over a single polygon, i.e. one point for each polygon edge [resp. for each polygon ray]. This, in general, leads to a further saving in computational cost, considering that any single edge integration point is shared by two adjacent polygons, hence is processed once on the overall computations. This of course implies to process the entire mesh in a modified manner, i.e. edge/polygon-ray-wise, instead of in an element-by-element fashion. The latter strategy has been adopted throughout the following VEM computations for $k=2$. 
We finally observe that even more efficient choices could be used by following the various results in the literature on integration rules for polygons, see for instance \cite{NME2759,Sommariva2007,Mousavi2011}.

\subsection{Thick-walled viscoelastic cylinder subjected to internal pressure}
\label{ss:viscocyl}
\begin{figure}
\begin{minipage}[b]{.5\linewidth}
\centering
    \includegraphics[bb=30 0 500 30, angle=0, scale=0.30]{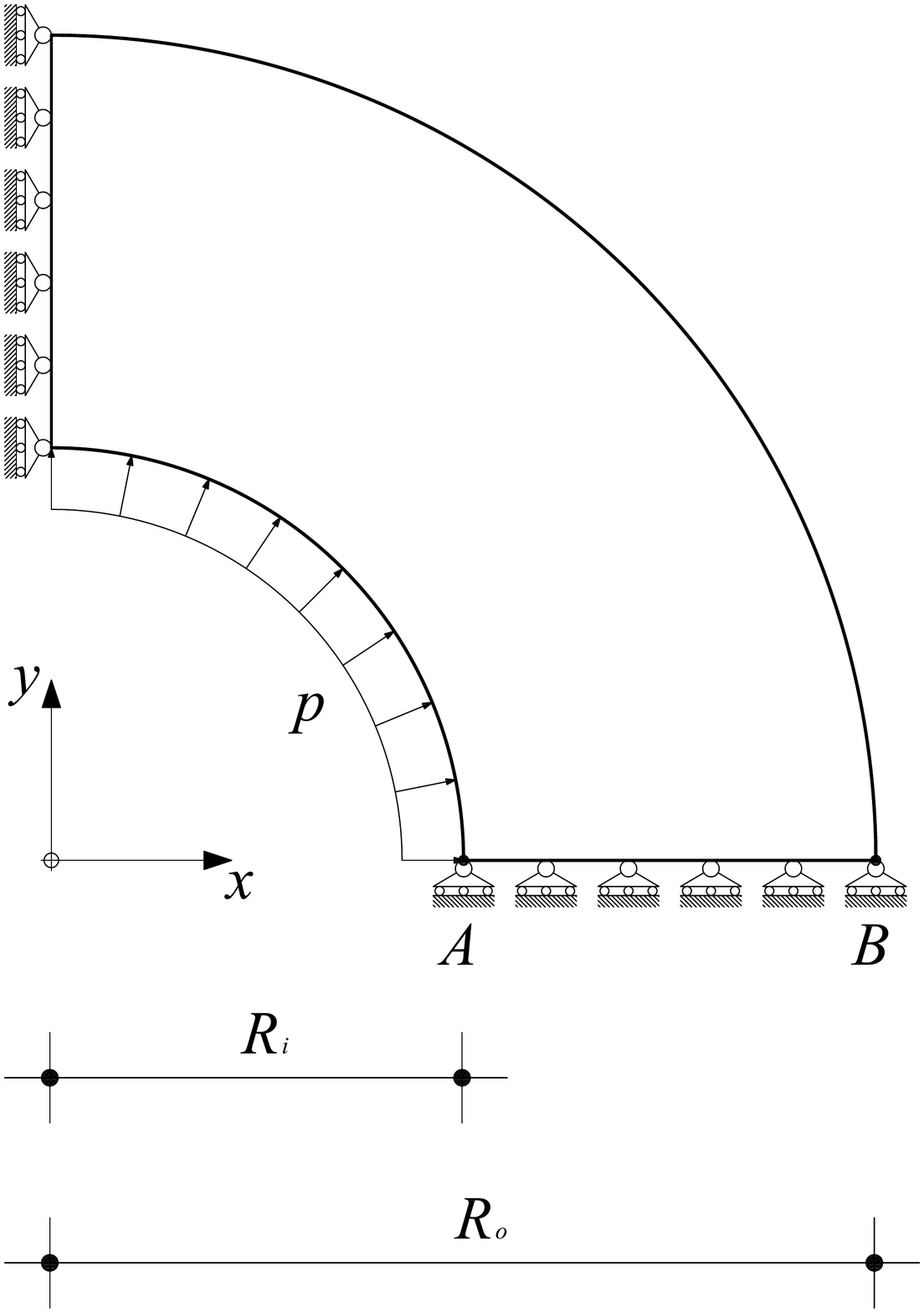}
\subcaption{}
\end{minipage}
\hspace{5mm}
\begin{minipage}[b]{.5\linewidth}
\centering
\includegraphics[bb=130 5 500 600, clip, angle=0, scale=0.425]{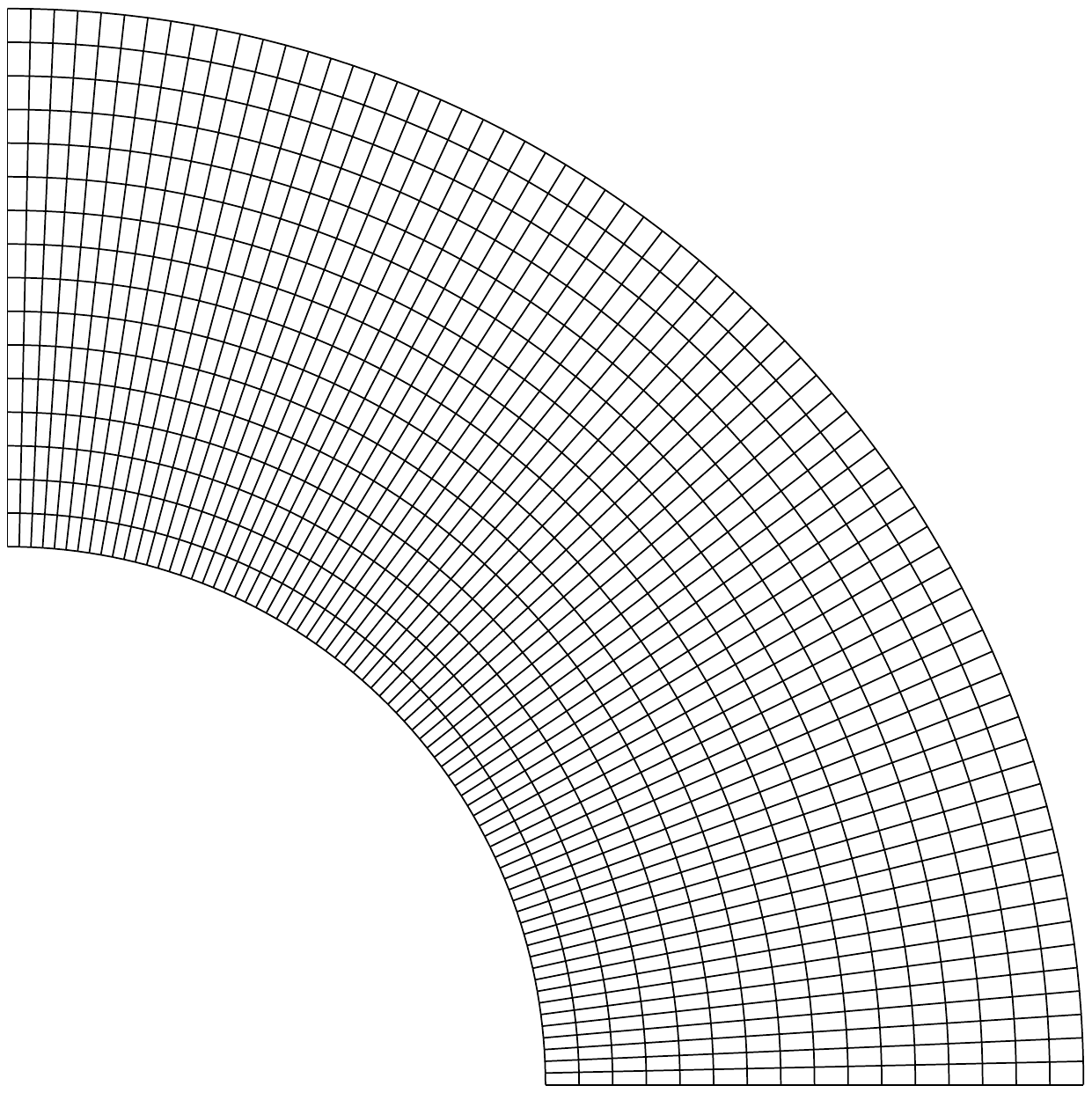}
\subcaption{}
\end{minipage}
\caption{Thick-walled viscoelastic cylinder subjected to internal pressure. (a) Geometry, boundary conditions, applied load; (b) structured quadrilateral mesh.}
\label{fig:pressure_cylinder_geom}
\end{figure}

The first benchmark regards a thick-walled cylinder characterized by a response reproduced by a viscoelastic constitutive model The cylinder, subjected to an internal pressure,  has inner [resp. outer] radius $R_i = 2$ [$R_o=4$]. For symmetry, only a quarter of the cylinder cross section is studied, as reported in Fig. \ref{fig:pressure_cylinder_geom} (a), imposing zero normal displacement along the radial edges. The material is considered to be isotropic and modeled by viscoelastic response in deviatoric stress-strain only, in compliance with the constitutive model outlined in Sec. \ref{ss:viscoelmodel}. The material properties are set assuming $M=1$ and $\lambda_1 \equiv  \lambda = 1$, i.e. a {\it standard linear solid} \cite{Zienkiewicz_Taylor_Zhu13} is considered. Young's modulus and Poisson's ratio are set $E = 1000$, $\nu = 0.3$, respectively. Two sets of viscoelastic parameters are adopted for the present analysis (cf. \eqref{eq:Gt}), i.e. $\left( \mu_0 , \mu_1 \right)_{\textrm{ve}1} = \left( 0.01, 0.99 \right )$ and, $\left (\mu_0 , \mu_1 \right)_{\textrm{ve}2} = \left( 0.3, 0.7 \right )$, respectively. 
The former case is calibrated in such a way that the ratio of the bulk modulus to shear modulus for instantaneous loading is given by $K/G(0) = 2.167$ and for long time loading, say at $t=8$, by $K/G(8) = 216.7$, which indicates a near incompressible behavior for sustained loading cases (at $t = \infty$ the Poisson ratio results $0.498$). The second material set indicates a sort of intermediate response for theoretically infinite time after loading application. 

The structural response for a suddenly applied internal pressure $p = 10$ is computed through $20$ unit time integrations. In particular, time integration of the constitutive equation is performed for a set of discrete points $t_k$, $k = 1,...,20$; by using the generalized Maxwell model in Prony series form (cf. Eq. \eqref{eq:Gt}), solution is reduced to a recursion formula in which each material state computed by a simple update of the previous one. Details concerning the implemented algorithm may be found in \cite{Zienkiewicz_Taylor_Fox}. 
It is noted that, albeit the adopted constitutive model is intrinsically three-dimensional, the above mentioned 2D solution procedure developed in the context of the plane strain assumption can carried out without any modification of the standard integration algorithm code illustrated in \cite{Zienkiewicz_Taylor_Fox}.

\begin{figure}
\begin{minipage}[b]{.5\linewidth}
\centering
\includegraphics[bb=90 150 500 600, clip, angle=0, scale=0.425]{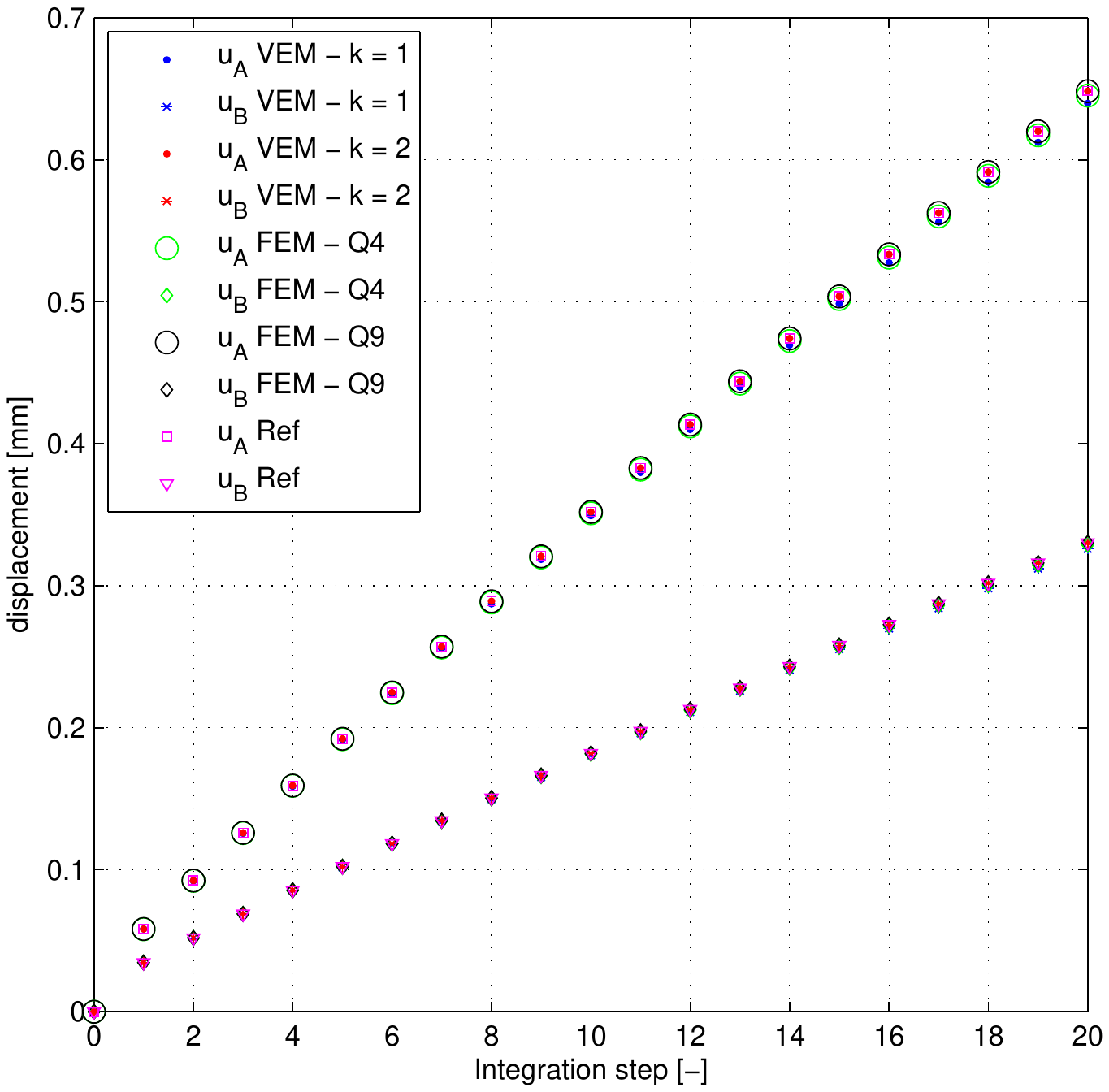}
\subcaption{}
\end{minipage}
\hspace{5mm}
\begin{minipage}[b]{.5\linewidth}
\centering
    \includegraphics[bb=90 150 500 600, clip, angle=0, scale=0.425]{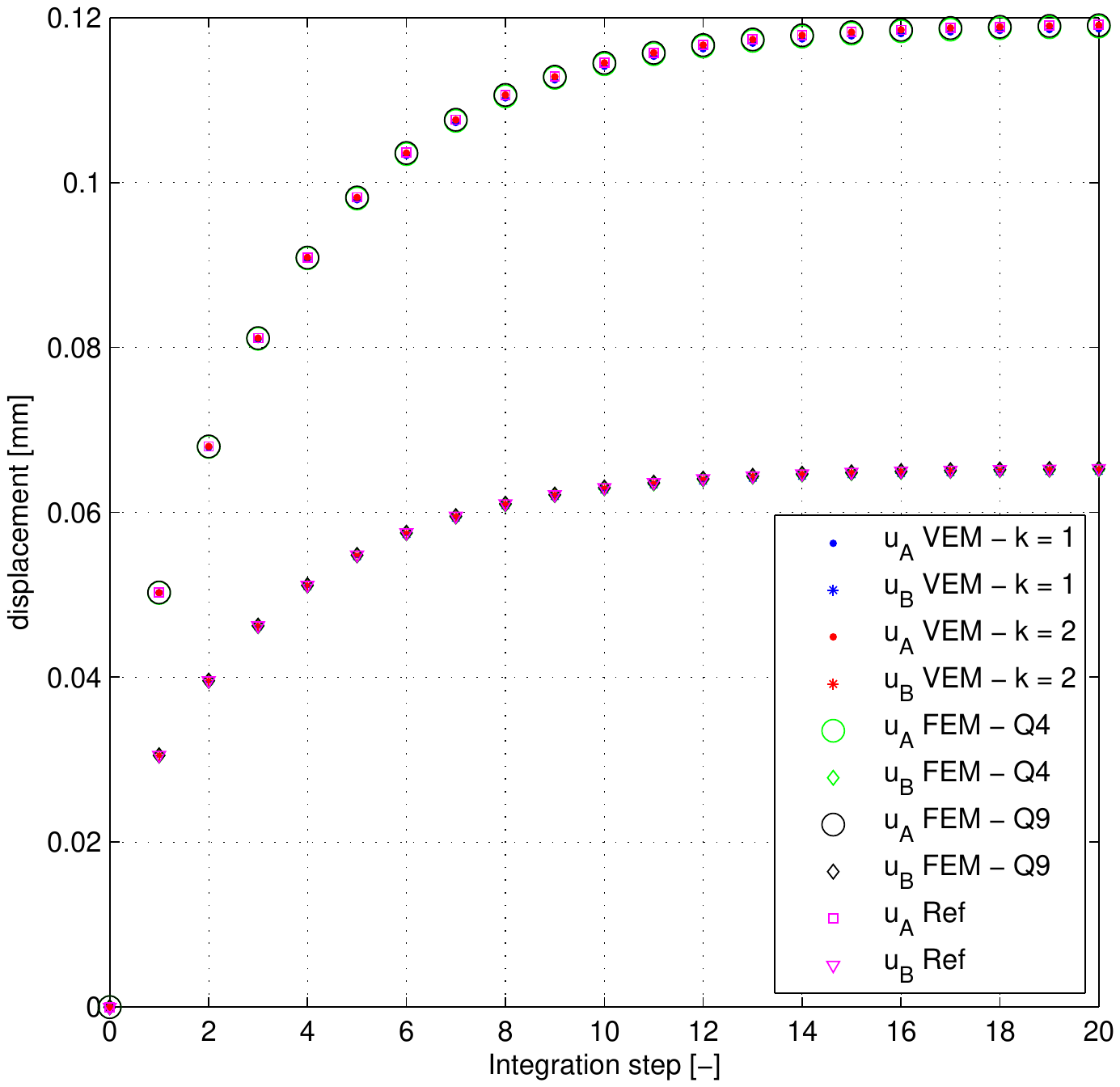}
\subcaption{}
\end{minipage}
\caption{Thick-walled viscoelastic cylinder with internal pressure. Integration step {\it vs.} displacement curve. (a) case $\left( \mu_0 , \mu_1 \right)_{\textrm{ve}1} = \left( 0.01, 0.99 \right )$; (b) case $\left (\mu_0 , \mu_1 \right)_{\textrm{ve}2} = \left( 0.3, 0.7 \right )$.}
\label{fig:pressure_cylinder_response}
\end{figure}

In the present work, a comparison is drawn between the proposed VEM formulation, for $k=1$ and $k=2$, and quadrilateral displacement based finite elements with four nodes $Q4$ (linear quadrilateral) and nine nodes $Q9$ (quadratic quadrilateral) \cite{Zienkiewicz_Taylor_Fox}. The adopted structured mesh is shown in Fig. \ref{fig:pressure_cylinder_geom} (b). A reference solution to this problem is computed with mixed $\bu-p-\beps^{\textrm{v}}$ quadrilateral finite elements and an overkilling space discretization \cite{Zienkiewicz_Taylor_Zhu13}. 
The integration-step versus the displacement curves for control points A and B (see Fig. \ref{fig:pressure_cylinder_geom} (a)) is shown in Fig. \ref{fig:pressure_cylinder_response} (a)-(b) for the compared solutions and for the reference one, for the two material parameter sets ${\textrm{ve}1}$ and ${\textrm{ve}2}$, introduced above. 
It is observed that the VEM formulation presents a response in excellent agreement with standard displacement finite elements for a given spatial discretization. 

This simple benchmark highlights indeed how implementation of the proposed VEM method into existing structural codes results straightforward and how numerical tools such as integration algorithms for viscoelastic constitutive equations translate immediately into the new framework without modification. 

\subsection{Perforated plastic plate}
\label{ss:plasticstrip}
\begin{figure}
\begin{minipage}[b]{.5\linewidth}
\centering
\includegraphics[bb=50 0 600 810, clip, angle=0, scale=0.3]{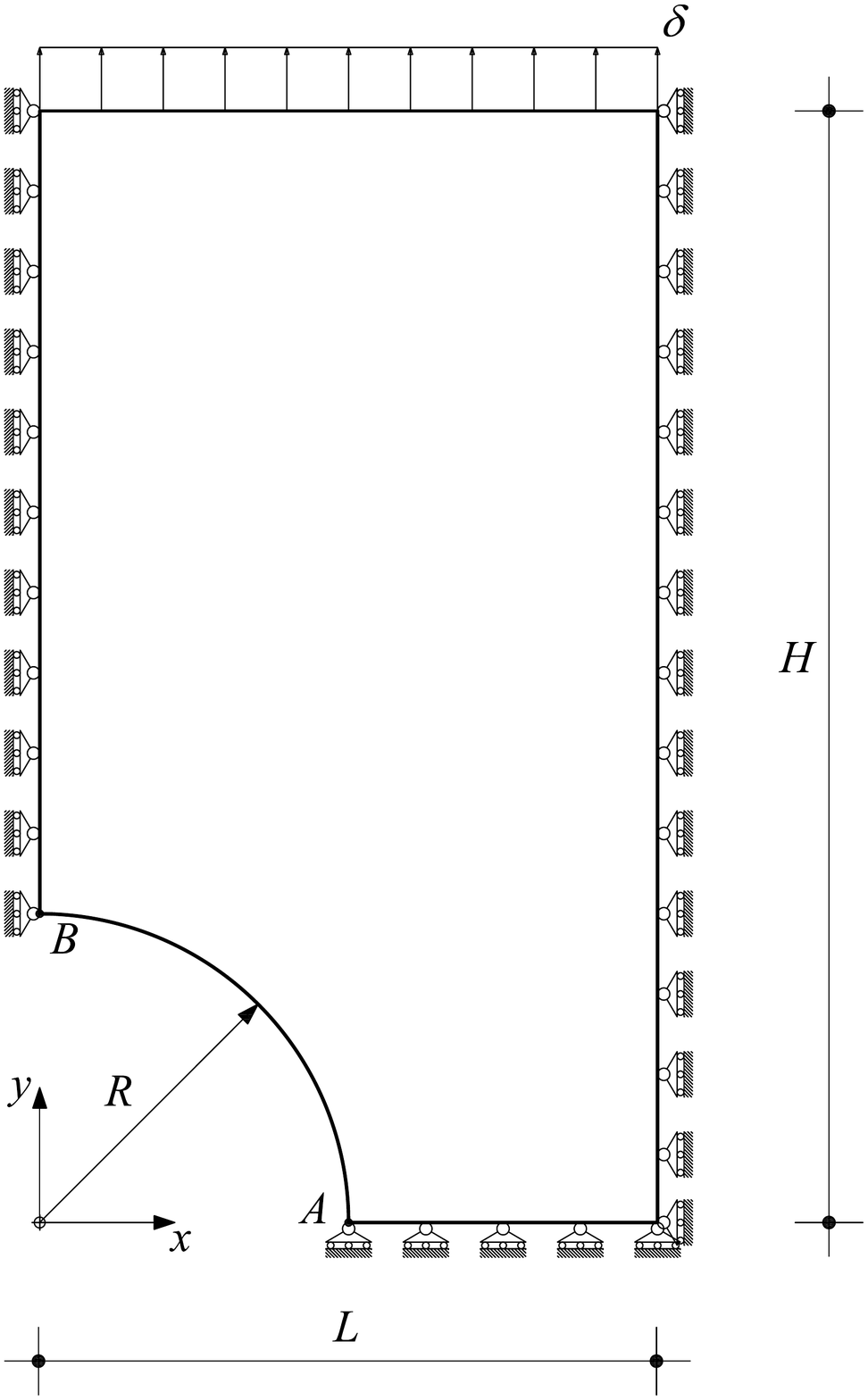}
\subcaption{}
\end{minipage}
\hspace{5mm}%
\begin{minipage}[b]{.5\linewidth}
\centering
\includegraphics[bb=130 155 450 600, clip, angle=0, scale=0.53]{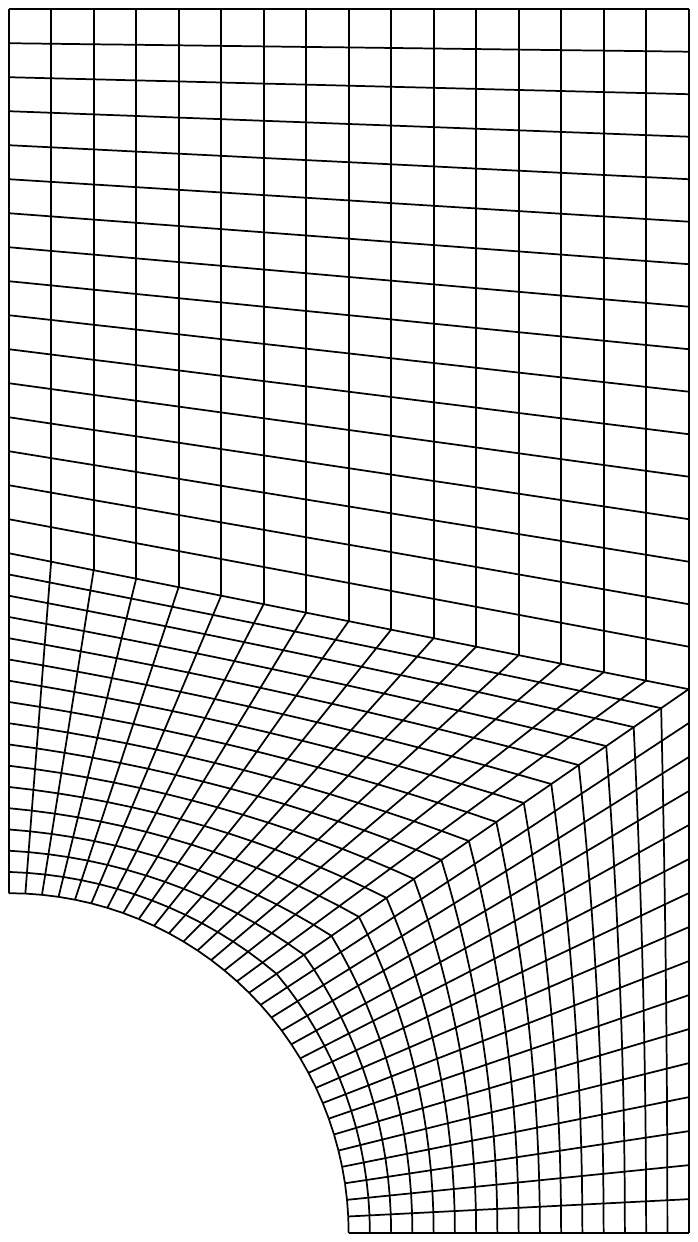}
\subcaption{}
\end{minipage}
\begin{minipage}[b]{.5\linewidth}
\centering
\includegraphics[bb=165 200 450 600, clip, angle=0, scale=0.53]{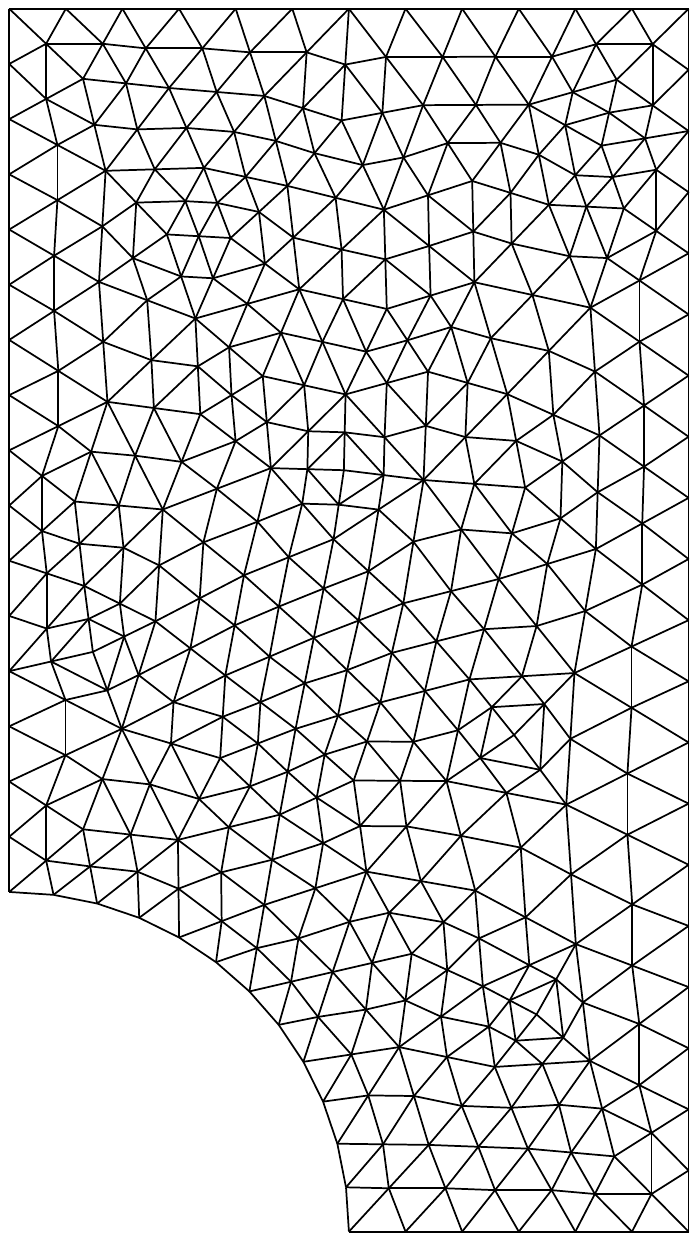}
\subcaption{}
\end{minipage}
\hspace{5mm}%
\begin{minipage}[b]{.5\linewidth}
\centering
\includegraphics[bb=130 200 450 600, clip, angle=0, scale=0.53]{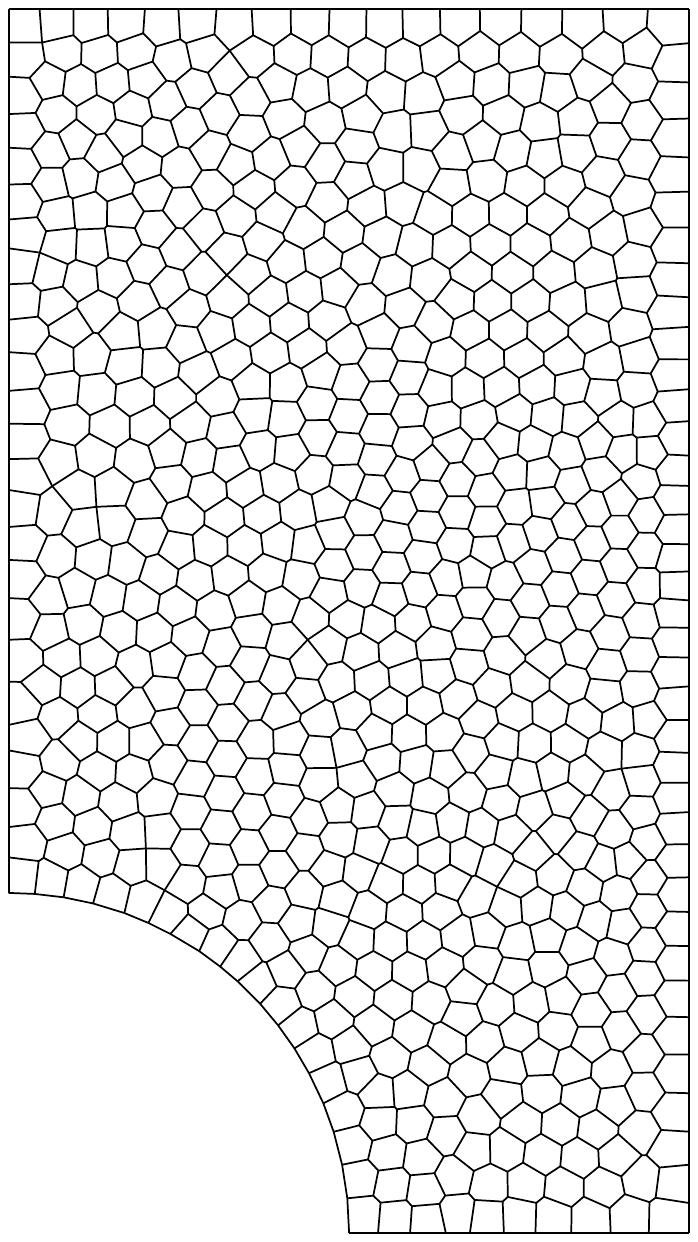}
\subcaption{}
\end{minipage}
\caption{Perforated plastic plate. (a) Geometry, boundary conditions, imposed displacement; (b) Quad - structured quadrilateral mesh; (c) Tri - unstructured triangular mesh; (d) Voronoi - centroid based tessellation.}
\label{fig:strip_wt_hole}
\end{figure}

A rectangular strip with $2L = 200$ $\textrm{mm}$ width and $2H = 360$ $\textrm{mm}$ length containing a central circular hole of $2R = 100$ $\textrm{mm}$ diameter is considered (viz. Fig. \ref{fig:strip_wt_hole} (a)). Material obeys to a plastic constitutive model (see Sec. \ref{ss:plastmodel}), and is characterized by $E = 7000$ $\textrm{kg} / \textrm{mm}m^2$, $\nu = 0.3$,  $\sigma_{\textrm{y},0} = 24.3$ $\textrm{kg} / \textrm{mm}^2$ \cite{Zienkiewicz_Taylor_Zhu13}. Plane strain assumption is here invoked as well, hence any (native 3D) integration algorithm for the constitutive model reviewed in Sec. \ref{ss:plastmodel} may be coherently utilized for the purpose of computing stress update at the integration point level in the VEM framework. To this end, for the present computation, a standard backward Euler scheme with return map projection is used \cite{simo_computational_1998}. 

Owing to symmetry, only one quadrant of the perforated strip is discretized as shown in Fig. \ref{fig:strip_wt_hole} (a). Displacement boundary restraints are prescribed for normal components on symmetry boundaries and on top and lateral boundaries. Loading is applied by a uniform normal displacement $\delta = 2$ $\textrm{mm}$ with $400$ equal increments on the upper edge, see Fig. \ref{fig:strip_wt_hole} (a). A quarter of the plate is meshed, respectively into quadrilateral (Quad), triangular (Tri), and polygonal (Voronoi) elements obtained with a centroid based tessellation (viz. Fig. \ref{fig:strip_wt_hole} (b)-(c)-(d)) generated by using the code \cite{TPPM12}. 
The simulation campaign is set considering virtual elements of order $k=1$ and $k=2$,  for the adopted meshes. For comparison and validation purposes, triangular [resp. quadrilateral] displacement based finite elements of type $T3$ (linear), $T6$ (quadratic) [resp. $Q4$ (linear), and $Q9$ (quadratic)] (\cite{Zienkiewicz_Taylor_Fox}) are used for the first two meshing layouts. A reference solution obtained with an overkilling discretization of $u-p-\beps^{\textrm{v}}$ quadrilateral mixed finite elements is computed for accuracy assessment. In passing, we note that VEM elements with $m=3$, $k=1$ and triangular Lagrangian finite elements $T3$ coincide, hence results pertaining only to the former scheme are reported in the following for conciseness.

\begin{figure}
\centering
\includegraphics[bb=50 150 550 810, clip, angle=0, scale=0.425]{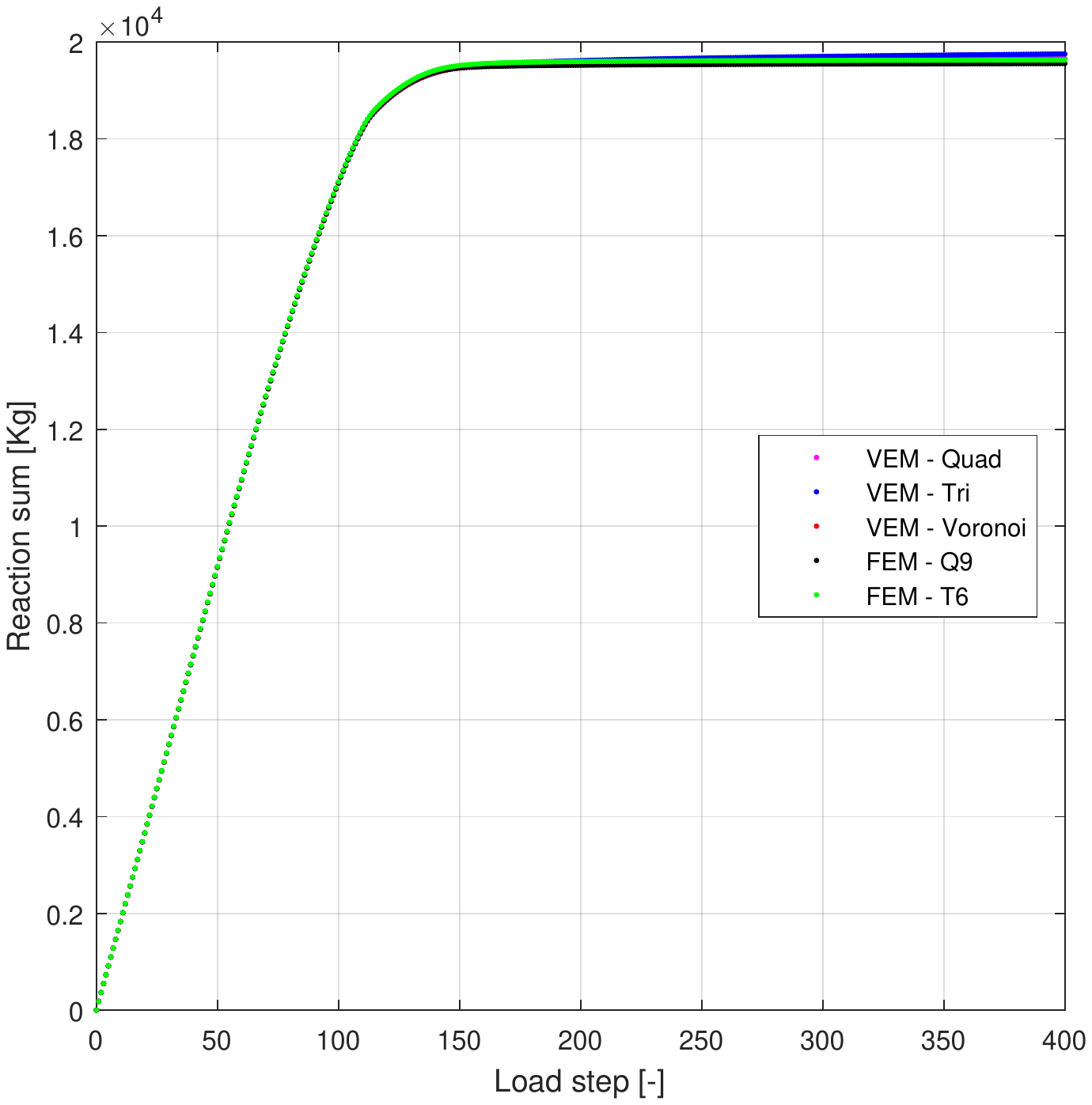}
\caption{Perforated plastic plate. Reaction sum {\it vs.} vertical displacement of upper edge curves. Comparison between VEM formulation and standard FEM for various meshes. Quadratic VEM ($k=2$) and quadratic Lagrangian FEM ($T6$, $Q9$) elements.}
\label{fig:strip_wt_hole_loadcurve}
\end{figure}

We report the horizontal [resp. vertical] displacement of point $A$ [resp. $B$] (see Fig. \ref{fig:strip_wt_hole} (a)) at the end of the loading history in Table \ref{tab:dispAB} for the compared different methods. A good agreement, for both components, obtained with linear and quadratic VEM/FEM formulations is observed by comparison with the reference mixed FEM solution. Table \ref{tab:convergence} shows the convergence characteristics of the proposed VEM method compared to the standard FEM in terms of average number of Newton iterations per incremental loading step, indicating a slight edge in terms of efficiency in favor of the VEM methodology. The overall response of the structure to the applied load is reported in Fig. \ref{fig:strip_wt_hole_loadcurve}, where load-displacement curves for the quadratic methods are reported, still in excellent agreement.

The present benchmark validates the proposed VEM methodology in terms of accuracy and efficiency with respect to standard FEM, indicating an outperformance of the former in terms of mesh versatility and convergence robustness. Implementability of the innovative VEM methodology in a standard nonlinear FEM analysis code is still retained.

\subsection{Shape memory alloy device}
\label{ss:smaactuator}
The present numerical test aims at proving the capability of the proposed VEM methodology in conjunction with complex highly nonlinear material behavior, as the one shown by shape memory alloy materials (cf. Sec. \ref{ss:SMA}). 

\begin{figure}
\begin{minipage}[b]{.5\linewidth}
\centering
\includegraphics[bb=30 215 580 610, clip, angle=0, scale=0.27]
                {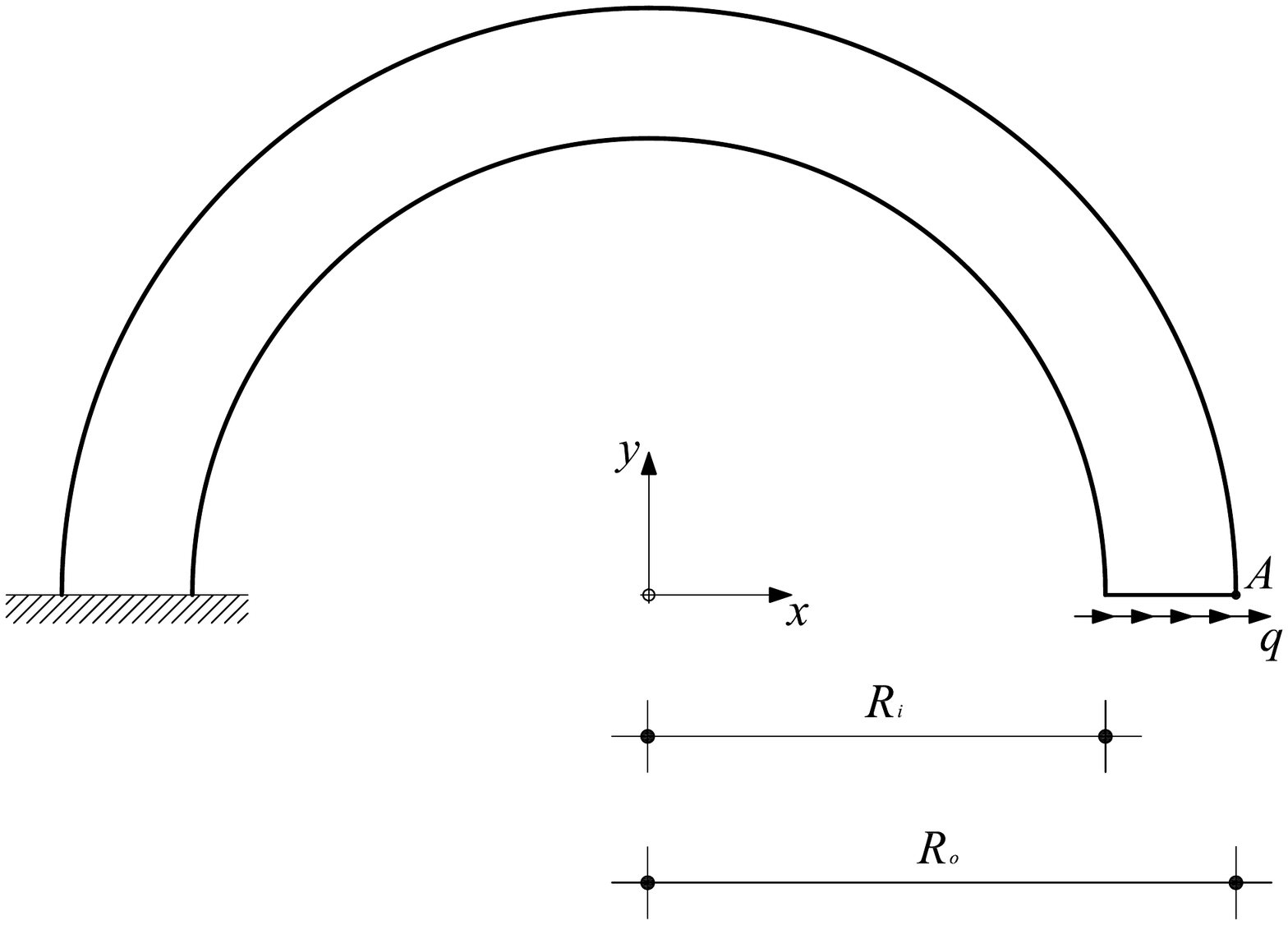}
\subcaption{}
\end{minipage}
\begin{minipage}[b]{.5\linewidth}
\centering
\includegraphics[bb=50 240 550 610, clip, angle=0, scale=0.34]
                {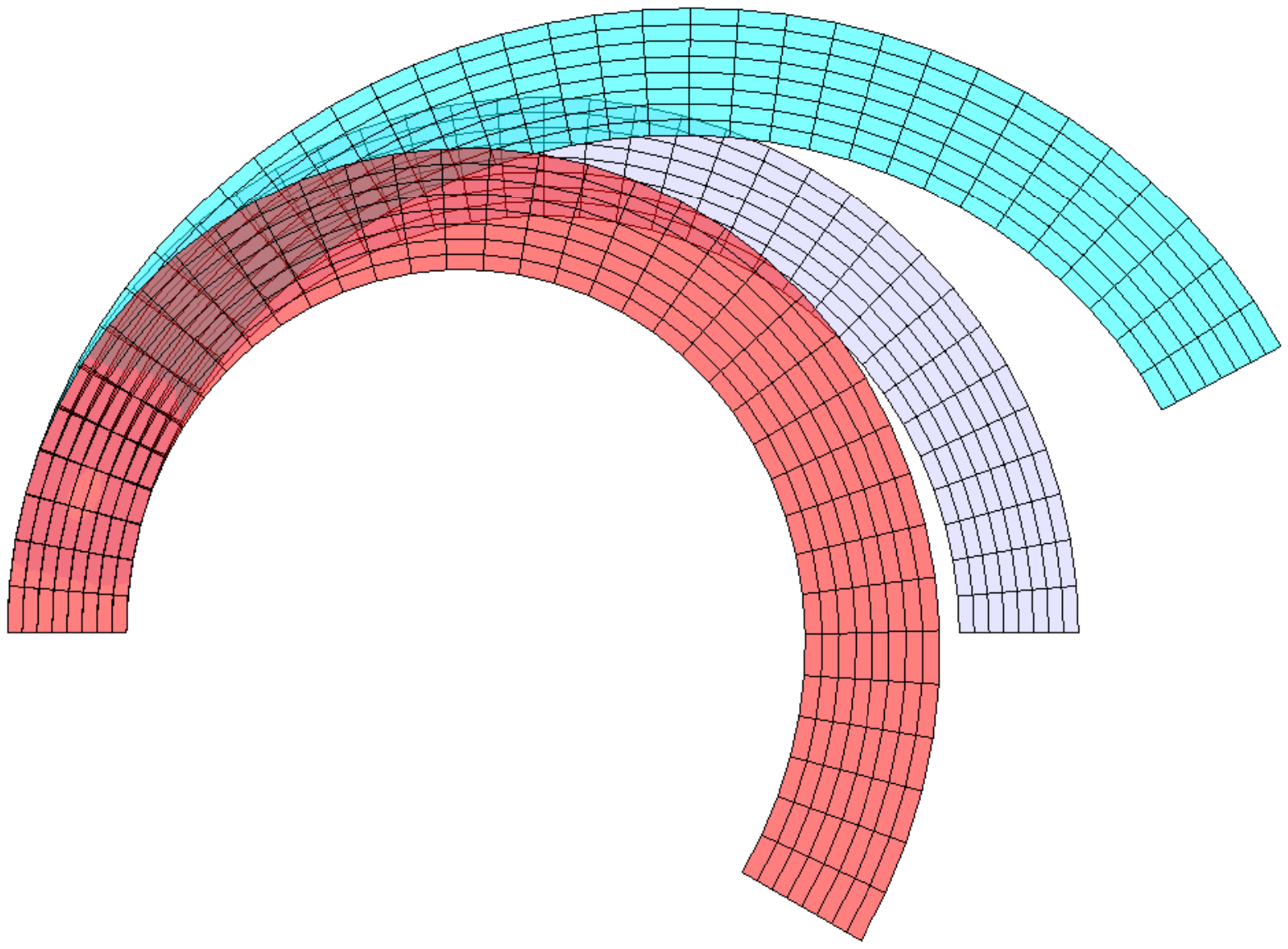}
\subcaption{}
\end{minipage}
\caption{Shape memory alloy device. (a) Geometry, boundary conditions, applied load; (b) adopted quadrilateral mesh: reference configuration (grey palette), deformed configurations at $t=1$ (cyan palette) and $t=3$ (red palette).}
\label{fig:SMA_actuator_geom_mesh}
\end{figure}

A typical clamped semicircular SMA arch device is considered (\cite{Evangelista2009,Evangelista2010,Artioli_Taylor_IJNME_2012}), as portrayed in Fig. \ref{fig:SMA_actuator_geom_mesh} (a) with indication of applied traction forces on the free end. Inner [resp. outer] radius is $R_{\textrm{i}} = 3.5$ $\textrm{mm}$ [resp. $R_{\textrm{o}} = 4.5$ $\textrm{mm}$]. Material parameters for the constitutive model are: $E = 53000$ $\textrm{MPa}$, $\nu = 0.36$, $\epsL = 0.04$, $\Mf = 223$ $\textrm{K}$, $h = 1000$ $\textrm{MPa}$, $\beta = 2.1$ $\textrm{MPa/K}$, $\sigma_{\textrm{y},0} = 50$ $\textrm{Mpa}$ \cite{Evangelista2009}. 

Differently from previous numerical tests, plane stress assumption is adopted in this case. The original 3D form of the constitutive model reviewed in Sec. \ref{ss:SMA} is integrated at the integration point level with the innovative state update algorithm recently proposed in \cite{Artioli_Bisegna_ECCM_2014,Artioli_Bisegna_IJNME_2016a}. To comply with plane stress condition, a further nonlinear constraint onto the stress state stemming from the state update is applied through a reduction algorithm with nested iterations, as outlined in \cite{Zienkiewicz_Taylor_Fox,SouzaNeto}, for instance. All in all, when looked at in a lower-to-upper operational level, the solution procedure amounts to an initial boundary problem on a 2D domain with triple nonlinearity, namely at: constitutive update, plane stress constraint prescription, and equilibrium problem solution level, respectively.

The structure is subjected to a load-temperature controlled loading history which comprises $5$ pseudo-time branches as indicated in table \ref{tab:SMAloadhist}, supposing proportional loading, with $40$ equal increments during the first $4$ branches and $N=10$ equal increments during the final heating branch. Loading parameters are $q_{\textrm{max}} = 60$ $\textrm{N} / \textrm{mm}$, $T_{\textrm{room}} = 223$ $\textrm{K}$. The VEM solution is obtained with the mesh of quad elements shown in Fig. \ref{fig:SMA_actuator_geom_mesh} (b) (grey palette). The results presented herein refer to the quadratic case, i.e. $k=2$. 

\begin{figure}
\centering
\includegraphics[bb=50 150 550 610, clip, angle=0, scale=0.425]
                {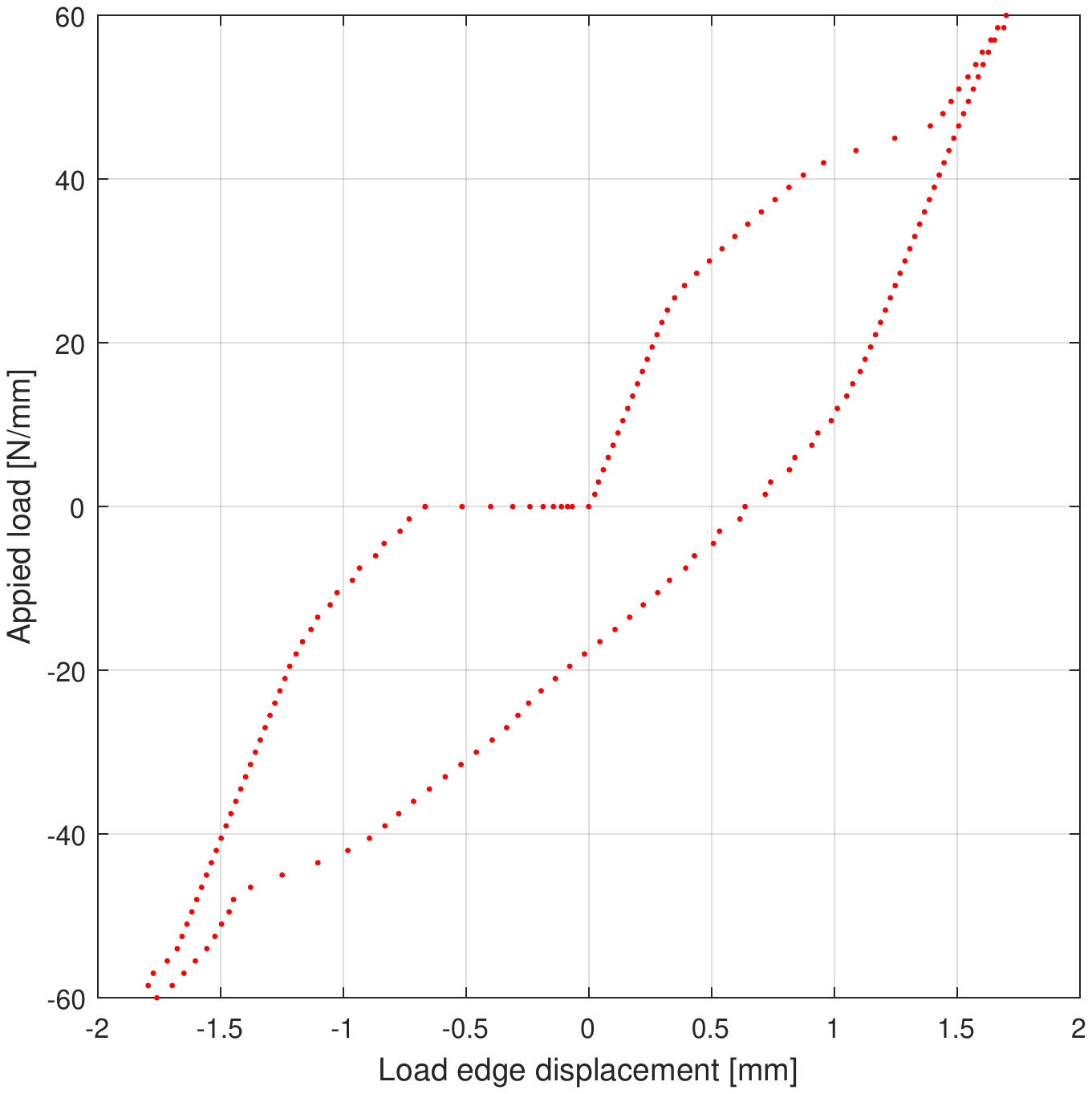}
\caption{Shape memory alloy device. Applied load {\it vs.} point $A$ displacement curve.}
\label{fig:SMA_actuator_curve}
\end{figure}

The deformed configurations on the adopted mesh at $t=1$ (cyan palette) and $t=3$ (red palette), respectively, are presented in Fig. \ref{fig:SMA_actuator_geom_mesh} (b). The functioning curve for the analyzed SMA device reports applied load {\it vs.} horizontal displacement of point A (cf. Fig. \ref{fig:SMA_actuator_geom_mesh} (a)) as can be seen in Fig. \ref{fig:SMA_actuator_curve}, where we can appreciate the classical shape recovery exhibited by the arch device. 

\section{Conclusion}
\label{s:conclusion}
The VEM formulation presented in Part I \cite{ABLS_part_I} for 2D elasticity problem, has been extended to the case of nonlinear material response. In particular, three classical and typical nonlinear material models have been considered, i.e. the Maxwell viscoplasticity, the Mises plasticity, and a SMA constitutive law.

It is remarked throughout the paper that the implementation of material nonlinearity laws are implemented in the VEM code substantially in the same way as in standard FEM framework. Thus, the solution algorithm for the typical integration point can be regarded as a black-box that can be extracted from classical FEM implementation and introduced in VEM code, without substantial changes.
In VEM formulation, the integration over domains characterized more than 3 edges has been performed dividing the polygon in triangles and using Gauss quadrature technique on each triangle.

Numerical results remark the ability of the VEM formulation to get accurate solutions for linear \cite{ABLS_part_I} and, now, also for nonlinear 2D structural problems. Solutions obtained using VEM are generally more accurate than the FEM ones and, important feature, require less iterations than FEM approach when the same tangent (Newton) algorithm is adopted.

Finally, the two parts of the present study, devoted to the development of a VEM formulation for linear and nonlinear 2D structural problems, demonstrate the accuracy of the approach and, also, the almost simplicity of the implementation, mainly when a FEM code is available. These two features make the VEM very interesting for a wide class of structural applications and, indeed, the possibility to the inclusion of VEM elements in available commercial FEM codes.

\clearpage
\newpage


\begin{table}[!htbp]
\centering
\caption{Perforated plastic plate. Convergence assessment in terms of displacement components $u_A$ and $v_B$ at the end of loading history.}
\begin{tabular}{lcccccc}
\toprule
& \multicolumn{5}{c}{linear elmts.}\\
\midrule
& VEM - Quad & VEM - Tri/FEM - T3 & VEM - Voronoi & FEM - Q4 & &Ref.\\
\midrule
$u_A$ $[\textrm{mm}]$ & 2.538  &  2.485  &  2.540  &  2.692 & &2.741 \\
$v_B$ $[\textrm{mm}]$ & 1.819  &  1.823  &  1.815  &  1.836 & &1.859 \\
\midrule
& \multicolumn{5}{c}{quadratic elmts.}\\
\midrule
& VEM - Quad & VEM - Tri & VEM - Voronoi & FEM - Q9 & FEM - T6 & Ref.\\
\midrule
$u_A$ $[\textrm{mm}]$ & 2.720  &  2.719  &  2.714  &  2.733 & 2.719 & 2.741 \\
$v_B$ $[\textrm{mm}]$ & 1.851  &  1.852  &  1.852  &  1.853 & 1.853 & 1.859 \\
\bottomrule
\end{tabular}
\label{tab:dispAB}
\end{table}

\begin{table}[!htbp]
\centering
\caption{Perforated plastic plate. Convergence comparison in terms of average number of iterations per load step.}
\begin{tabular}{lccccc}
\toprule
& VEM - Quad & VEM - Tri & VEM - Voronoi & FEM - Q$\sharp$ & FEM - T$\sharp$\\
\midrule
\textrm{linear elmts.} & 4.86  &  5.04  &  4.63  &  5.44 & 5.04\\
\textrm{quadratic elmts.} & 6.43  &  6.11  &  6.21  &  6.70 & 6.11\\
\bottomrule
\end{tabular}
\label{tab:convergence}
\end{table}

\clearpage
\newpage

\begin{table}[!htbp]
\centering
\caption{Shape memory alloy device. Loading history chart for applied load and temperature.}
\begin{tabular}{l c ll ll ll l}
\toprule
Time $[-]$ & &$0$ & $1$ & $2$ & $3$ & $4$ & $5$\\
\midrule
\textrm{load} $q$ $[\textrm{N} / \textrm{mm}]$ & & $0$  & $q^\textrm{max}$  &  0  &  $-q^\textrm{max}$  &  0  & 0\\
\textrm{Temperature} $T$ $[\textrm{K}]$ & &$T^\textrm{room}$  &  $T^\textrm{room}$  &  $T^\textrm{room}$  &   $T^\textrm{room}$ &   $T^\textrm{room}$ & $T^\textrm{room} + 80$ \\
\bottomrule
\end{tabular}
\label{tab:SMAloadhist}
\end{table}

\clearpage
\newpage
%

\clearpage
\newpage

\clearpage
\newpage

\clearpage
\newpage

\clearpage
\newpage

\clearpage
\newpage

\bibliographystyle{spmpsci}      
\bibliography{VEM_Bibliography,SMA_Bibliography,general-bibliography,biblio}   

\end{document}